\documentclass[11pt]{article}
\usepackage{amsmath}
\usepackage{amsfonts}
\usepackage{amssymb}
\usepackage{graphicx}
\usepackage{color}

\numberwithin{equation}{section}

\definecolor{mycolorred}{rgb}{1, 0, 0}

\newtheorem{theorem}{Theorem}[section]

\newtheorem{corollary}[theorem]{Corollary}

\newtheorem{definition}[theorem]{Definition}

\newtheorem{lemma}[theorem]{Lemma}

\newtheorem{proposition}[theorem]{Proposition}
\newtheorem{remark}[theorem]{Remark}

\def\<{\langle}
\def\>{\rangle}

\def\E{{\mathbb E}}
\def\R{{\mathbb R}}

\def\N{{\mathbb N}}

\begin{document}

\title{An Invariance Principle for Stochastic Series II.\\
Non Gaussian Limits}
\author{ \textsc{Vlad Bally}\thanks{%
Universit\'e Paris-Est, LAMA (UMR CNRS, UPEMLV, UPEC), INRIA, F-77454
Marne-la-Vall\'ee, France. Email: \texttt{bally@univ-mlv.fr}.} \smallskip \\
\textsc{Lucia Caramellino}\thanks{%
Dipartimento di Matematica, Universit\`a di Roma - Tor Vergata, Via della
Ricerca Scientifica 1, I-00133 Roma, Italy. Email: \texttt{%
caramell@mat.uniroma2.it}}\smallskip\\
}
\maketitle

\parindent 0pt


\parindent 0pt

\begin{abstract}
We study the convergence in total variation
distance for series of the form
\begin{equation*}
S_{N}(c,Z)=\sum_{l=1}^{N}%
\sum_{i_{1}<...<i_{l}}c(i_{1},...,i_{l})Z_{i_{1}}...Z_{i_{l}}
\end{equation*}
where $Z_{k},k\in {\mathbb{N}}$ are independent centered random variables
with ${\mathbb{E}}(Z_{k}^{2})=1. $ This enters in the framework of the $U$--%
statistics theory which plays a major role in modern statistic. In the case
when $Z_{k},k\in {\mathbb{N}}$ are standard normal, $S_{N}(c,Z)$ is an
element of the sum of the first $N$ Wiener chaoses and, starting with the
seminal paper of D. Nualart and G. Peccati, 
the convergence of
such functionals to the Gaussian law has been extensively studied. So the
interesting point consists in studying invariance principles, that is, to
replace Gaussian random variables with random variables with a general law.
This has been done in several papers using the Fortet--Mourier distance, the
Kolmogorov distance or the total variance distance. In particular, estimates
of the total variance distance in terms of the fourth order cumulants has
been given in the part I of the present paper. 
But, as the celebrated Fourth Moment Theorem of
Nualart and Peccati 
shows, such estimates are pertinent to deal with Gaussian limits.
In the present paper we study the convergence to
general limits which may be non Gaussian, and then the estimates of the
error has to be done in terms of the low influence factor only.
\end{abstract}

\smallskip


\tableofcontents

\section{Introduction and main results}

Let us introduce the objects involved in our paper. We consider a sequence
of independent random variables $Z_{k},k\in {\mathbb{N}}$ with ${\mathbb{E}}%
(Z_{k})=0$ and ${\mathbb{E}}(Z_{k}^{2})=1.$ We assume that the law of each
of them is locally lower bounded by the Lebesgue measure, that is ${\mathbb{P%
}}(Z_{k}\in dy)\geq \varepsilon dy$ for $y\in B(z_{k},2r).$ More precisely,
there exists $r,\varepsilon >0$ and $z_{k}\in {\mathbb{R}}$ such that, for
every measurable function $f:{\mathbb{R}}\rightarrow {\mathbb{R}}_{+}$%
\begin{equation}
{\mathbb{E}}(f(Z_{k}))\geq \varepsilon \int f(z)1_{B(0,2r)}(z-z_{k})dz.
\label{Int2}
\end{equation}

All along the paper we will fix some $r,\varepsilon \in (0,1)$ and an
increasing sequence $M_{p}\in (1,\infty ),p\in {\mathbb{N}}.$ These are
arbitrary but fixed (without any supplementary mention). We use the notation $%
\mathcal{L}((M_{p})_{p\in {\mathbb{N}}},r,\varepsilon )$ to indicate the sequences of
independent random variables $Z=(Z_{k})_{k\in {\mathbb{N}}}$\ with ${\mathbb{%
E}}(Z_{k})=0$ and ${\mathbb{E}}(Z_{k}^{2})=1$ which verifies (\ref{Int2})
with $r,\varepsilon $ and such that $\left\Vert Z_{k}\right\Vert _{p}\leq
M_{p}$ for every $k,p\in {\mathbb{N}}.$ Notice that the random variables $%
Z_{k}$ are not identically distributed. However, the fact that we may choose
$(M_{p})_{p\in {\mathbb{N}}},r,\varepsilon $ to be the same for all of them
represents an uniformity property.

We consider a family of coefficients $c=\{c(\alpha ):\alpha \in {\mathbb{N}}%
^{m},m\in {\mathbb{N}}\}$ and for a multi-index $\alpha =(\alpha
_{1},...,\alpha _{m})\in\N^m$ we denote $\left\vert \alpha \right\vert =m$ the
length of $\alpha .$ We also denote $Z^{\alpha }=Z_{\alpha _{1}}\cdots Z_{\alpha _{m}}.$ We denote by $\mathcal{C}$ the class of the
coefficients $c$ which are symmetric and null on the diagonals. And we look
to stochastic series of the following type:%
\begin{equation}
S_{N}(c,Z)=\sum_{1\leq \left\vert \alpha \right\vert \leq N}c(\alpha
)Z^{\alpha }  \label{Int1}
\end{equation}%
This enters in the framework of $U$--statistics introduced by Hoeffding \cite%
{Hoe} and Fisher \cite{Fish}, which play a major role in modern statistics
(see for example Lee \cite{lee}). Moreover we denote%
\begin{equation}\label{Int3}
\delta _{1}(c) =\max_{k}\left\vert c(k)\right\vert ,\quad \delta
_{m}(c)=\max_{k}\Big(\sum_{\left\vert \alpha \right\vert =m-1}c^{2}(\alpha
,k)\Big)^{1/2}\  m\geq 2,\quad
\overline{\delta }_{N}(c) =\sum_{m=1}^{N}\delta _{m}(c).
\end{equation}
$\overline{\delta }_{N}(c)$ is the so called ``influence factor'': $\sum_{\left\vert \alpha
\right\vert =m-1}c^{2}(\alpha ,k)$ may be considered as the measure of the
action of the particle $k$ on all the other particles, at level $m.$ And if $%
\overline{\delta }_{N}(c)$ is small we say that we have ``low influence''.

We will also use the following semi-norms%
\begin{equation}
\left\vert c\right\vert _{m}=\Big(\sum_{\left\vert \alpha \right\vert
=m}c^{2}(\alpha )\Big)^{1/2}\quad \mbox{and}\quad \left\Vert c\right\Vert
_{N}^{2}=\sum_{m=1}^{N}\left\vert c\right\vert _{m}^{2}  \label{not11}
\end{equation}

We are now able to give our first result:

\begin{theorem}
\label{CONVERGENCE} We consider a sequence $Z^{n}=(Z_{k}^{n})_{k\in {\mathbb{%
N}}}\in \mathcal{L}((M_{p})_{p\in {\mathbb{N}}},r,\varepsilon ).$ Let $N\in {%
\mathbb{N}}$ be fixed an let $c_{n}\in \mathcal{C}$ be a sequence of
coefficients such that%
\begin{equation}
\limsup_{n}\sum_{\left\vert \alpha \right\vert \leq N}c_{n}^{2}(\alpha
)<\infty .  \label{Int4}
\end{equation}%
We assume that they verify the ``low influence condition'':%
\begin{equation}
\lim_{n\rightarrow \infty }\overline{\delta }_{N}(c_{n})=0.  \label{Int5}
\end{equation}%
We also assume that the following non degeneracy condition holds:
\begin{equation}
\liminf_{n\rightarrow \infty }\sum_{\left\vert \alpha \right\vert
=N}c_{n}^{2}(\alpha )>0.  \label{ND}
\end{equation}%
Let $\mu $ be a probability measure. Then $\lim_{n\rightarrow \infty
}S_{N}(c_{n},Z)=X $ in law implies (and so is equivalent to) convergence
in total variation distance.
\end{theorem}

\begin{remark}
This is a generalization of the celebrated Prohorov's Theorem (see \cite%
{proh}) concerning convergence in total variation in the CLT (which
corresponds to $N=1).$ And as it is clear from Prohorov's theorem, the
condition (\ref{Int2}) appears as natural when dealing with convergence in
total variation distance (in contrast with convergence in law or in
Kolmogorov distance when such a condition is not necessary). A more
particular variant of this result has already been obtained recently by
Nourdin and Poly in \cite{[NPy]}.
\end{remark}

\begin{remark}
Notice that the non degeneracy condition (\ref{ND}) is much stronger than
the one in \cite{[Numelin1]} where $\sum_{\left\vert \alpha \right\vert
=N}c_{n}^{2}(\alpha )$ is replaced by $\sum_{\left\vert \alpha \right\vert
\leq N}c_{n}^{2}(\alpha ).$ So here we ask that the higher line of $S_{N}$
is non degenerated while in \cite{[Numelin1]} all the coefficients $c(\alpha
)$ in the sum contribute to the non degeneracy condition. But there we also
need that the cumulants tend to zero (not only the influence factor) and if
this is true then $\mu $ is a Gaussian probability measure.
\end{remark}

We will now give some (non asymptotic) estimates for the errors involved in
the limit in total variation distance. We denote ${\mathbb{N}}_{\ast}$ the
set of the positive integers and given $N\in {\mathbb{N}}_{\ast }$ we will
use the following constants:%
\begin{equation}
c_{N}(r,\varepsilon )=\Big(\frac{\varepsilon \sqrt r}{\sqrt{2}}\Big)^{2N}%
\frac{1}{N}  \label{C1}
\end{equation}%
and we use the generic notation $C_{N}(r,\varepsilon )$ for every constant
of the form
\begin{equation}
C_{N}(r,\varepsilon )=C(N!)^{q_{1}}e^{q_{2}M_{p}}\,r^{-q_{3}}\varepsilon
^{-q_{2}}  \label{C2}
\end{equation}%
where $C,p,q_{i}\in {\mathbb{N}}_{\ast },i=1,...,4$ are universal constants
(independent of the parameters $M_{p},\varepsilon ,r$ and on $N)$ and which
may change from a line to another.

We first estimate the error which is done by replacing a sequence $%
Z=(Z_{k})_{k\in {\mathbb{N}}}$ with another sequence $\overline{Z}=(%
\overline{Z}_{k})_{k\in {\mathbb{N}}}:$ this is the invariance principle. We
recall first Theorem 3.1
from \cite{[Numelin1]} which which concerns smooth
test functions (notice that here the hypothesis (\ref{Int2}) is not
necessary):

\begin{theorem}
\label{Smooth}Let $Z=(Z_{k})_{k\in {\mathbb{N}}}$ and $\overline{Z}=(%
\overline{Z}_{k})_{k\in {\mathbb{N}}}$ be two sequences of centered
independent random variables such that ${\mathbb{E}}(Z_{k}^{2})={\mathbb{E}}(%
\overline{Z}_{k}^{2})=1.$ We also assume that ${\mathbb{E}}(\left\vert
Z_{k}\right\vert ^{3})\leq M_{3}$ and ${\mathbb{E}}(\left\vert \overline{Z}%
_{k}\right\vert ^{3})\leq M_{3}.$ Then for every $f\in C_{b}^{3}({\mathbb{R}}%
)$ and every $c\in \mathcal{C}$%
\begin{equation}
\left\vert {\mathbb{E}}(f(S_{N}(c,Z)))-{\mathbb{E}}(f(S_{N}(c,\overline{Z}%
)))\right\vert \leq (N+1)!^{2}M_{3}^{4N}\left\Vert f^{\prime \prime \prime
}\right\Vert _{\infty }\|c\|_N\overline{\delta }_{N}(c).  \label{C3}
\end{equation}
\end{theorem}

The aim of the present paper is to obtain a similar estimate but to replace $%
\left\Vert f^{\prime \prime \prime }\right\Vert _{\infty }$ by $\left\Vert
f\right\Vert _{\infty },$ that is to work in total variation distance. This
has already been done in \cite{[Numelin1]} (see Theorem 6.1 therein) but
there the estimate involves the fourth cumulant (and not only $\overline{%
\delta }_{N}(c)).$ So, if we aim to use such estimates in order to study the
convergence of a sequence $S_{N}(c_{n},Z),n\in {\mathbb{N}},$ then the limit
has to be a Gaussian random variable (this is a consequence of the Fourth
Moment Theorem of Nualart and Peccati \cite{[NP]}). In the present paper we
prove the following estimate in terms of $\overline{\delta }_{N}(c)$ (which
is allows to study the convergence to general laws):

\begin{theorem}
\label{INVARIANCE}
Let $Z=(Z_{k})_{k\in {\mathbb{N}}}$ and $\overline{Z}=(%
\overline{Z}_{k})_{k\in {\mathbb{N}}}$ be two sequences of random variables
which belong to $\mathcal{L}((M_{p})_{p\in {\mathbb{N}}},r,\varepsilon )$
and let $c\in \mathcal{C}$. Then, for every $N$ and for every bounded and
measurable function $f:{\mathbb{R}}\rightarrow {\mathbb{R}}$%
\begin{equation}\label{Int6}
\begin{array}{rl}
\left\vert {\mathbb{E}}(f(S_{N}(c,Z)))-{\mathbb{E}}(f(S_{N}(c,\overline{Z}%
)))\right\vert \leq &C_{N+1}(r,\varepsilon )\left\Vert f\right\Vert
_{\infty } (1+\|c\|_N)\smallskip \\
&\displaystyle
\times \Big(\frac{\overline{\delta }_{N}^{\frac{1}{4+3p_\ast N}}(c)}{|c|_N^{\frac{6p_\ast}{4+3p_\ast N}}}
+\exp \Big(-\frac{%
c_{N}(r,\varepsilon )\left\vert c\right\vert _{N}^{2}}{\overline{\delta }%
_{N}^{2}(c)}\Big)\Big),
\end{array}%
\end{equation}
with $\left\Vert c\right\Vert _{N}$ and $\left\vert c\right\vert _{N}$
defined in (\ref{not11}) and $p_{\ast }$ is the universal constant which
appears in (\ref{not14}).
\end{theorem}

Similar but less precise results have been obtained before. Assume for a
moment that we replace $S_{N}(c,Z)$ by $\Phi _{N}(c,Z):=\sum_{\left\vert
\alpha \right\vert =N}c(\alpha )Z^{\alpha }.$ A first result, concerning
convergence in law, has been obtained in the pioneering papers of de Jong
\cite{[dJ1],[dJ2]}. Afterwards, in \cite{[MO'DO]} the authors prove
convergence in Kolmogorov distance, that is
\begin{equation*}
\sup_{x}\left\vert {\mathbb{E}}(1_{(-\infty ,x)}(\Phi _{N}(c,Z)))-{\mathbb{E}%
}(1_{(-\infty ,x)}(\Phi _{N}(c,\overline{Z})))\right\vert \rightarrow 0\quad
\mbox{as}\quad \delta _{N}(c)\rightarrow 0.
\end{equation*}%
These results hold for general random variables $Z_{k},$ condition (\ref%
{Int2}) being not needed. And recently, Nourdin and Poly in \cite{[NPy]}
assume (\ref{Int2}), and they prove that
\begin{equation*}
\sup_{\left\Vert f\right\Vert _{\infty }\leq 1}\left\vert {\mathbb{E}}%
(f(\Phi _{N}(c,Z)))-{\mathbb{E}}(f(\Phi _{N}(c,\overline{Z})))\right\vert
\rightarrow 0\quad \mbox{as}\quad \delta _{N}(c)\rightarrow 0.
\end{equation*}

A first progress in our paper is that we consider a general sum $S_{N}(c,Z)$
and not only $\Phi _{N}(c,Z).$ And more important, we obtain an estimate of
the error - and this is not asymptotic, but holds for every fixed $c\in
\mathcal{C}.$

The drawback of the estimate (\ref{Int6}) is that it rapidly degradates as $N
$ becomes large. This point is a consequence of the techniques we use here: we
use a stochastic variation calculus (analogues to the Malliavin calculus)
and the delicate point is to estimate the Malliavin covariance matrix
associated to our series; in order to do this we use Carbery-Wright
inequality which concerns general polynomials and which make appear $1/N$ as
a power of $\overline{\delta }_{N}(c).$ One may compare this estimate with
the similar one which is given in Theorem 6.2
in \cite{[Numelin1]}. There the upper bound is given in terms of the fourth cumulant $%
\kappa _{4,N}(c)$ of $\Phi _{N}(c,G),$ where $G=(G_{k})_{k\in {\mathbb{N}}}$
with $G_{k}$ are independent standard normal random variables. And that
upper bound is of the form $C_{N}(r,\varepsilon )\kappa _{4,N}^{1/2}(c)$ when $c(\alpha)=0$ for $|\alpha|=1$, otherwise the power is no more $1/2$ but $1/4$.
In any case, the power of $\kappa _{4,N}(c)$ does not depend on $N$. However we stress
that the two estimates may not be directly compared because $\overline{%
\delta }_{N}(c)\leq \kappa _{4,N}(c),$ and it is possible that $\overline{%
\delta }_{N}(c)$ is much smaller than $\kappa _{4,N}(c)$ (see e.g. the example developed in in Section \ref{ex2}).

The estimate of the Malliavin covariance matrix is done in \cite{[Numelin1]}
using some martingale techniques which take into account the specific
structure of the stochastic series at hand and so are more powerful than
estimates concerning general polynomials (as in the Carbery-Wright
inequality). But they make appear the fourth cumulant $\kappa _{4,N}(c)$
which does not converge to zero, except in the case when we focus on a
Gaussian limit (as it is pointed out by the fourth moment theorem of Nualart
and Peccati \cite{[PN]}). So, if we aim to general limits, we have to come
back to the Carbery-Wright lemma (which does not involve cumulants).

We give now some estimates of the error in the convergence in total
variation of a sequence $S_{N}(c_{n},Z),n\in {\mathbb{N}}$ to a probability
measure $\mu .$ We will work with the metrics
\begin{equation}
d_{k}(F,G)=\sup \{\left\vert {\mathbb{E}}(f(F))-{\mathbb{E}}%
(f(G))\right\vert :\left\Vert f\right\Vert _{k,\infty }\leq 1\}  \label{Int7}
\end{equation}%
where
\begin{equation}
\Vert f\Vert _{k,\infty }:=\sum_{p=0}^{k}\Vert
f^{(p)}\Vert _{\infty }.  \label{Int8}
\end{equation}%
In particular $d_{0}=d_{TV}$ is the total variation distance and $%
d_{1}=d_{FM}$ is the Fortet Mourier distance (which metrizes the convergence
in law).

\begin{theorem}
\label{DISTANCES-2}
Let $X$ be a random variable which is the limit in law of  $(S_M(c_n,\overline{Z}))_n$, for some $M\in \N_*$, where $(\overline{Z}_{k})_{k\in \N}\in \mathcal{L}((M_{p})_{p\in{\mathbb{N}}},r,\varepsilon )$ and, for every $n\in\N$,  $c_{n}\in\mathcal{C}$ satisfies (\ref{Int4}), (\ref{Int5}) and (\ref{ND}) with $N$ replaced by $M$. Set
\begin{equation}\label{CNX}
\overline{C}_{M,X}=\limsup_n\|c_n\|_M\quad\mbox{and}\quad \underline{C}_{M,X}=\liminf_n|c_n|_M.
\end{equation}
Then for every $c\in\mathcal{C}$ and $(Z_{k})_{k\in {\mathbb{N}}}\in \mathcal{L}((M_{p})_{p\in\N},r,\varepsilon )$ one has
\begin{equation} \label{Int9-2}
\begin{array}{rl}
d_{0}(S_N(c,Z),X ) \leq &
\displaystyle
C_{N\vee M}(r,\varepsilon )(1+\| c\|_{N}+\overline{C}_{M,X})\\
&\displaystyle
\Big(
\frac{d_{1}^{\frac 1{2+p_{\ast }N\vee M}}(S_N(c,Z),X) }
{(|c|_N^{2/N}\wedge \underline{C}_{M,X}^{2/M})^{\frac{p_\ast N\vee M}{2+p_{\ast }N\vee M}}}
+\exp \Big(-\frac{c_{N}(r,\varepsilon )\left\vert c\right\vert _{N}^{2}}{%
\delta _{N}^{2}(c)}\Big)\Big),
\end{array}
\end{equation}
$c_N(r,\varepsilon)$ and $C_N(r,\varepsilon)$ being as in (\ref{C1}) and (\ref{C2}) respectively and $p_{\ast }$ is the uni\-ver\-sal constant from (\ref{not14}).
\end{theorem}

We can rewrite Theorem \ref{DISTANCES-2} by using the concept of ``$M$--attainability''.

\begin{definition}
\label{Ateignable}Given $M\in {\mathbb{N}}$ and $(M_{p})_{p\in\N},\varepsilon
,r>0$ we say that $X$ is $M$-attainable of class $(M_{p})_{p\in\N},%
\varepsilon ,r$ if there exists a sequence of coefficients $c_{n}\in
\mathcal{C}$ which satisfy (\ref{Int4}),(\ref{Int5}) and (\ref{ND}) with $N$ replaced by $M$ and a
sequence $Z^{n}=(Z_{k}^{n})_{k\in {\mathbb{N}}}\in \mathcal{L}((M_{p})_{p\in
{\mathbb{N}}},r,\varepsilon )$ such that $\lim_{n}S_{M}(c_{n},Z^{n})=X $
in law. If $X$ is $M$--attainable, we set $\overline{C}_{M,X}$ and $\underline{C}_{M,X}$ as in (\ref{CNX}).

We denote by $\mathcal{A}_{M}((M_{p})_{p\in\N},\varepsilon ,r)$ this
class.
\end{definition}

If $M=1$ the CLT for non identically distributed random variables shows that
the only $1$-attainable random variable is the standard normal one. And if $%
M=2$ a characterization of the $2$-attainable laws is given in \cite{[NPy2]}
(see also \cite{Sev}). Of course they include random variables equal in law to elements in the second chaos. And more generally, elements in a fixed $M$ chaos are $M$-attainable.
So, as an immediate consequence of Theorem \ref{DISTANCES-2} we obtain the following

\begin{corollary}
\label{DISTANCES}Let $Z\in \mathcal{L}((M_{p})_{p\in {\mathbb{N}}%
},r,\varepsilon ),$ $c\in \mathcal{C}$ and $X \in \mathcal{A}%
_{M}((M_{p})_{p\geq 1},r,\varepsilon )$. Let $p_{\ast }$
be the uni\-ver\-sal constant from (\ref{not14}). Then
$d_{0}(S_M(c,Z),X ) $ satisfies inequality (\ref{Int9-2}).
\end{corollary}

The proofs of the above results are given in Section \ref{proofs}.

Finally we give several examples of applications.

First, in Theorem \ref{CHI2}, we estimate the distance between $\Phi
_{N}(c,Z)$ and a $\chi _{2}$ law with $m$ degrees of freedom. This
significantly straighten a result of Nourdin and Peccati from \cite{[NP2009]}
concerning approximation of the law of a multiple stochastic integral by a $%
\chi _{2}$ law with $m$ degrees of freedom: the result in \cite{[NP2009]} concerns
Wiener multiple integrals and the estimate is in $d_{1}$ distance, while
here we have a general sequence of random variables $Z_{k}$ and the estimate
is in terms of $d_{0}.$

In a second application we prove that
\begin{equation*}
S_{n}(Z)=\frac{n!}{\sqrt{2n\ln n}}\sum_{1<i<j\leq n}\frac{1}{\sqrt{j-i}}%
Z_{i}Z_{j}
\end{equation*}%
converges to the standard normal distribution and the total variation
distance to the limit is upper bounded by $(n^{-1}\ln
^{2}n)^{1/43(1+2p_{\ast })}.$ We notice that if the interaction potential $%
|j-i|^{-1/2}$ is replaced by $|j-i|^{-p}$ with $p\in (0,1/2)$, then (with a
suitable renormalization) the above sum converges to a double stochastic
integral. So, if $p=1/2$ we have a contraction phenomenon whereas such a
phenomenon does not exist if $p<1/2$.

Finally, in the third example we consider
\begin{equation*}
X_{i}=\frac{1}{\sqrt{n}}\sum_{\substack{ j=1  \\ j\neq i}}^{n}\frac{1}{\sqrt{%
\left\vert j-i\right\vert }}Z_{j}\quad \mbox{and}\quad
V_{n}(Z)=\sum_{i=1}^{n}X_{i}^{2}-EX_{i}^{2}
\end{equation*}%
and we prove that $V_{n}(Z)$ converges to a double Wiener integral and the
total variation distance is upper bounded by $(n^{-1}\ln
^{2}n)^{1/4(3+2p_{\ast })}.$

Finally in Appendix \ref{hoef} we prove an iterated version of Hoeffding's
inequality (which may be of own interest) and in Appendix \ref{computations} we give
some estimates for integrals needed in the last two examples.

\section{Proofs of the main results}

\subsection{Notation and preliminary results}

All along we consider some sequence $(M_{p})_{p\in {\mathbb{N}}}$ and some $%
\varepsilon ,r>0$ to be given, and we employ the notation already settled in
the introduction. We consider a sequence of random variables $%
Z=(Z_{k})_{k\in {\mathbb{N}}}\in \mathcal{L}((M_{p})_{p\geq 1},r,\varepsilon
),$ so each $Z_{k}$ satisfies (\ref{Int2}). Then we construct a function $%
\psi _{r}$ in the following way:%
\begin{equation}
\theta _{r}(z)=1-\frac{1}{1-(\frac{z}{r}-1)^{2}}\qquad \psi
_{r}(z)=1_{\{\left\vert z\right\vert \leq r\}}+1_{\{r<\left\vert
z\right\vert \leq 2r\}}e^{\theta _{r}(\left\vert z\right\vert )}.
\label{not2}
\end{equation}%
We denote%
\begin{equation}
m(r)=\int_{{\mathbb{R}}}\psi _{r}(\left\vert z\right\vert ^{2})dz\leq 2\sqrt{%
2r},\qquad v(r)=\frac{1}{m(r)}\int_{{\mathbb{R}}}z^{2}\psi _{r}(\left\vert
z\right\vert ^{2})dz\geq \frac{r}{3\sqrt 2}.  \label{not3}
\end{equation}%
Then $m(r)^{-1}\psi _{r}(\left\vert z\right\vert ^{2})$ is a probability
density and the corresponding random variable has mean zero and variance $%
v(r).$

Since $\psi _{r}\leq 1_{B(0,2r)},$ the inequality (\ref{Int2}) holds with $%
1_{B(0,2r)}$ replaced by $\psi _{r}.$ This allows to use a splitting method
in order to give the following representation of the law of $Z_{k}.$ We
consider some independent random variables $\chi _{k},U_{k},V_{k}$, $k\in{%
\mathbb{N}}$, with
\begin{eqnarray*}
{\mathbb{P}}(\chi _{k} &=&1)=\varepsilon m(r),\qquad {\mathbb{P}}(\chi
_{k}=0)=1-\varepsilon m(r), \\
{\mathbb{P}}(U_{k} &\in &dz)=\frac{1}{m(r)}\psi _{r}(\left\vert
z-z_{k}\right\vert ^{2})dz \\
{\mathbb{P}}(V_{k} &\in &dz)=\frac{1}{1-\varepsilon m(r)}({\mathbb{P}}%
(Z_{k}\in dz)-\varepsilon \psi _{r}(\left\vert z-z_{k}\right\vert ^{2})dz.
\end{eqnarray*}%
(we stress that in in \cite{[Numelin1]} the role of $U$ and $V$ are inverted).

Then $\chi _{k}U_{k}+(1-\chi _{k})V_{k}$ has the same law as $Z_{k}$ so from
now on we assume that
\begin{equation*}
Z_{k}=\chi _{k}U_{k}+(1-\chi _{k})V_{k}.
\end{equation*}

We will work with stochastic series based on $Z_{k},$ which we introduce
now. We denote $\Gamma _{m}={\mathbb{N}}_{\ast }^{m}$. Any $\alpha\in\Gamma_m
$ is named a multi-index and we define $\left\vert \alpha \right\vert =m$
its length. We set $\Gamma=\cup_{m}\Gamma_m$. For $J\in {\mathbb{N}}_*$ we
denote $\Gamma _{m}(J)=\{\alpha \in \Gamma _{m}:1\leq \alpha _{i}\leq J\}$
and $\Gamma (J)=\cup _{m=1}^{\infty }\Gamma _{m}(J).$ Moreover, for $%
z_{i}\in {\mathbb{R}},i\in {\mathbb{N}}$ and $\alpha =(\alpha
_{1},...,\alpha _{m})\in\Gamma_m$, we denote $z^{\alpha }=\prod_{i=1}^{m}z_{\alpha
_{i}}.$ We denote by $\mathcal{C}$ the class of the coefficients $%
c\,:\,\Gamma\to {\mathbb{R}}$ which are symmetric and null on the diagonals.
Then we consider a family of coefficients $c\in \mathcal{C}$ and we work
with the stochastic series%
\begin{equation}
S_{N}(c,Z)=\sum_{1\leq \left\vert \alpha \right\vert \leq N}c(\alpha
)Z^{\alpha }=\sum_{m=1}^{N}\sum_{\alpha \in \Gamma _{m}}c(\alpha )Z^{\alpha
}.  \label{not4}
\end{equation}%
In \cite{[Numelin1]} we developed a stochastic variational calculus based on
$U_{k},k\in {\mathbb{N}}$ (the explicit expression of the density of the law
of $U_{k}$ is central in that calculus) but here we do not need to recall
all this -- we will just recall some consequences which are used in the
present paper. We denote%
\begin{equation}
c_{j}(\alpha )=(1+\left\vert \alpha \right\vert )c(\alpha ,j)  \label{not5}
\end{equation}%
with the convention that, if $\alpha $ is void, then $\left\vert \alpha
\right\vert =0$ and $c_{j}(\alpha )=c(j).$ Then we define
\begin{equation}
\lambda _{N}=\lambda _{S_{N}(c,Z)}:=\sum_{j=1}^{\infty }\chi _{j}\left\vert
\partial _{Z_{j}}S_{N}(c,Z)\right\vert ^{2}=\sum_{j=1}^{\infty }\chi
_{j}\left\vert c(j)+S_{N-1}(c_{j},Z)\right\vert ^{2}.  \label{not6}
\end{equation}%
This is the ``Malliavin covariance matrix'' (in our one-dimensional case,
this is a scalar) associated to $S_{N}(c,Z)$ and plays a central role in our
estimates. Moreover we recall the seminorms $|c|_m$ and $\|c\|_N$ in (\ref%
{not11}) and we define
\begin{equation}
N_{q}(c,M)=\Big(\sum_{m=q}^{N}M^{m-q}\times \frac{m!}{(m-q)!}\times
m!\sum_{\left\vert \alpha \right\vert =m}c^{2}(\alpha )\Big)^{1/2}\leq N!e^{%
\frac{1}{2}M}\left\Vert c\right\Vert _{N}.  \label{not10}
\end{equation}%
where $\left\Vert c\right\Vert _{N}$ is defined in (\ref{not11}). In \cite%
{[Numelin1]} (see (4.17) therein
we have defined the
Sobolev norms $\left\Vert \left\vert S_{N}(c,Z)\right\vert \right\Vert _{q,p}
$ and in \cite{[Numelin1]} Proposition 5.5,
formula (5.14),
we have proved that
\begin{equation*}
\left\Vert \left\vert S_{N}(c,Z)\right\vert \right\Vert _{q,p}\leq \frac{C}{%
r^{q-1}}\Big(\sum_{n=1}^{q-1}\left\vert c\right\vert _{n}+N_{q}(c,M_{p}^{2})%
\Big).
\end{equation*}%
So, using (\ref{not10}) we have
\begin{equation}
\left\Vert \left\vert S_{N}(c,Z)\right\vert \right\Vert _{q,p}\leq \frac{C}{%
r^{q-1}}\,N!e^{M_p/2}\,\|c\|_N.  \label{not12}
\end{equation}%
In fact the only way in which $\left\Vert \left\vert S_{N}(c,Z)\right\vert
\right\Vert _{q,p}$ comes on in the present paper is just by means of the
above inequality, so the reader does not need to go further in the knowledge
of this quantity.

We use now a regularization lemma from \cite{[Numelin1]}. Let $\psi _{1}$ be
the function defined in (\ref{not2}), $m(1)$ the normalization constant from
(\ref{not3}) (with $r=1)$ and, for $\delta >0,$ let%
\begin{equation*}
\gamma _{\delta }(z)=\frac{1}{m(1)\sqrt{\delta }}\psi _{1}(\delta
^{-1}\left\vert z\right\vert ^{2}).
\end{equation*}
For $f\,:\,\R\to\R$ we set $f\ast \gamma_\delta$ the convolution between $f$ and $\gamma_\delta$, whenever it is well defined.
Using the regularization Lemma 4.6
from \cite{[Numelin1]} and (\ref{not12}) we obtain

\begin{lemma}
\label{L3}There exist some universal constants $C,p_{\ast }\geq 1$ such that
for every $\eta >0,\delta >0$ and for every bounded and measurable $f:{%
\mathbb{{\mathbb{R}}}}\rightarrow {\mathbb{R}}$ one has%
\begin{equation}
\left\vert {\mathbb{E}}(f(S_{N}(c,Z)))-{\mathbb{E}}(f\ast \gamma _{\delta
}(S_{N}(c,Z)))\right\vert \leq C_{N}(r,\varepsilon )\left\Vert c\right\Vert
_{N}\left\Vert f\right\Vert _{\infty }\Big({\mathbb{P}}(\lambda _{N}<\eta )+%
\frac{\sqrt{\delta }}{\eta ^{p_{\ast }}}\Big)  \label{not14}
\end{equation}%
with $C_{N}(r,\varepsilon )=Cr^{-2}N!e^{\frac 12M_{p_{\ast }}}.$
\end{lemma}

We will use the following easy consequence, which is a slightly more precise
version of Theorem 2.7 from \cite{[BC-EJP]}.

\begin{lemma}
\label{d0dk} Let $Z,\overline{Z}\in \mathcal{L}((M_{p})_{p\geq
1},r,\varepsilon )$ and $c,\overline{c}\in \mathcal{C}.$\ Let $p_{\ast }$ be
the universal constant from (\ref{not14}). For every $k\in {\mathbb{N}}$, $N,M\in\N_*$
there exists a universal constant $C_{N\vee M}(r,\varepsilon )$ (depending on $k)$
such that for every $\eta >0$ one has
\begin{equation}
\begin{array}{rcl}
d_{0}(S_{N}(c,Z),S_{M}(\overline{c},\overline{Z})) & \leq & \displaystyle %
C_{N\vee M}(r,\varepsilon )\big(1+\left\Vert c\right\Vert _{N}+\left\Vert \overline{c%
}\right\Vert _{M}\big)\times\smallskip \\
&  & \times \displaystyle \Big(\frac{1}{\eta ^{\frac{kp_{\ast}}{k+1}}}d_{k}^{\frac{1}{k+1%
}}(S_{N}(c,Z),S_{M}(\overline{c},\overline{Z})) +{\mathbb{P}}(\lambda
_{N}<\eta )+{\mathbb{P}}(\overline{\lambda }_{M}<\eta )\Big)%
\end{array}
\label{not16}
\end{equation}%
where $\lambda _{N}=\lambda _{S_{N}(c,Z)}$ and $\overline{\lambda }%
_{M}=\lambda _{S_{M}(\overline{c},\overline{Z})}$ are defined in (\ref{not6}%
), $d_{k}$ is defined in (\ref{Int7}) and $C_N(r,\varepsilon)$ is a constant
of the form (\ref{C2}).
\end{lemma}

\textbf{Proof.} To simplify the notation we put $S_N=S_{N}(c,Z)$ and $%
\overline{S}_M=S_{M}(\overline{c},\overline{Z})$ and $C$ a constant of the
form $C_{N\vee M}(r,\varepsilon )(1+\left\Vert c\right\Vert _{N}+\left\Vert
\overline{c}\right\Vert _{M})$ (which changes from a line to another). Let
$\delta >0$ and let $f\in C({\mathbb{R}})$ with $\left\Vert f\right\Vert
_{\infty }\leq 1.$ Since $\left\Vert f\ast \gamma _{\delta }\right\Vert
_{k,\infty }\leq C\delta ^{-k/2}$ we have
\begin{equation*}
\left\vert {\mathbb{E}}(f\ast \gamma _{\delta }(S_N))-{\mathbb{E}}(f\ast
\gamma _{\delta }(\overline{S}_M))\right\vert \leq C\delta ^{-k/2}d_{k}(S_N,%
\overline{S}_M).
\end{equation*}%
Then, using (\ref{not14}),%
\begin{equation*}
\left\vert {\mathbb{E}}(f(S_N))-{\mathbb{E}}(f(\overline{S}_M))\right\vert \leq
C\delta ^{-k/2}d_{k}(S_N,\overline{S}_M)+C\Big({\mathbb{P}}(\lambda _{N}<\eta )+
{\mathbb{P}}(\overline{\lambda }_{M}<\eta )+\frac{\delta ^{1/2}}{\eta
^{p_{\ast }}}\Big).
\end{equation*}%
We optimize over $\delta $: we take%
\begin{equation*}
\delta ^{(k+1)/2}=d_{k}(S,\overline{S})\eta ^{p_{\ast }}.
\end{equation*}%
We insert this in the previous inequality and we obtain (\ref{not16}). $%
\square $

\subsection{Estimate of the covariance matrix}

Our aim is to estimate ${\mathbb{P}}(\lambda _{N}\leq \eta )$ with $\lambda
_{N}$ defined in (\ref{not6}) and this will be done using the Carbery-Wright
inequality (we follow here an idea from \cite{[NPy]}). In order to
do this we need the following lemma.

\begin{lemma}
We denote by ${\mathbb{E}}_{V,\chi }$ the conditional expectation with
respect to $\sigma (V_{i},\chi _{i},i\in {\mathbb{N}}).$ Let $v(r)$ be as in
(\ref{not3}). Then%
\begin{equation}
{\mathbb{E}}_{V,\chi }(\lambda _{N})\geq \frac{v^{N-1}(r)}{N}%
\sum_{\left\vert \alpha \right\vert =N}c^{2}(\alpha )\chi ^{\alpha }.
\label{NDC6}
\end{equation}
\end{lemma}

\textbf{Proof.} We denote%
\begin{equation*}
\overline{U}_{i}=U_{i}-{\mathbb{E}}(U_{i})\qquad \mbox{and}\qquad \overline{V%
}_{i}=(1-\chi _{i})V_{i}+\chi _{i}{\mathbb{E}}(U_{i})
\end{equation*}%
so that%
\begin{equation*}
Z_{i}=\chi _{i}U_{i}+(1-\chi _{i})V_{i}=\chi _{i}\overline{U}_{i}+\overline{V%
}_{i}.
\end{equation*}%
Then we define%
\begin{equation*}
\overline{Z}^{\alpha }=\sum_{\substack{ (\beta ,\gamma )=\alpha ,  \\ \gamma
\neq \emptyset }}\chi ^{\beta }\overline{U}^{\beta }\times \overline{V}%
^{\gamma }
\end{equation*}%
and we write%
\begin{equation*}
Z^{\alpha }=\overline{Z}^{\alpha }+\chi ^{\alpha }\overline{U}^{\alpha }.
\end{equation*}%
Notice that for every multi-indexes $\alpha \in \Gamma _{m}$ with $m\leq N-1$
and $\theta \in \Gamma _{N}$ we have%
\begin{equation}
{\mathbb{E}}_{V,\chi }(\overline{Z}^{\alpha }\overline{U}^{\theta })=\sum
_{\substack{ (\beta ,\gamma )=\alpha ,  \\ \gamma \neq \emptyset }}\chi
^{\beta }{\mathbb{E}}_{V,\chi }(\overline{U}^{\beta }\overline{U}^{\theta
})\times \overline{V}^{\gamma }=0.  \label{NDCW7}
\end{equation}%
This is because $\left\vert \beta \right\vert <m\leq N-1$, so there is at
least one $\theta _{i}\notin \beta $ and ${\mathbb{E}}_{V,\chi }(\overline{U}%
^{\theta _{i}})=0.$ We take now $\kappa\in{\mathbb{R}}$ and we consider the
r.v. $X=\kappa+S_N(c,Z)$. We write write $\kappa+S_{N}(c,Z)=S^{\prime
}+S^{\prime \prime }$ with $S^{\prime }=\sum_{\left\vert \alpha \right\vert
=N}c(\alpha )\chi ^{\alpha }\overline{U}^{\alpha }$ and $S^{\prime \prime
}=\kappa+S_{N}(c,Z)-S^{\prime }.$ By (\ref{NDCW7}), $S^{\prime }$ and $%
S^{\prime \prime }$ are orthogonal in $L^{2}({\mathbb{P}}_{V,\chi })$ so
that
\begin{equation}
{\mathbb{E}}_{V,\chi }((\kappa+S_{N}(c,Z))^{2})\geq {\mathbb{E}}_{V,\chi }({%
S^{\prime }}^{2})=\sum_{\left\vert \alpha \right\vert
=N}v^{N}(r)c^{2}(\alpha )\chi ^{\alpha },  \label{NDCW7''}
\end{equation}%
the last equality being a consequence of ${\mathbb{E}}(\overline{U}%
_{i}^{2})=v(r).$

Consider now $\lambda_N$: by (\ref{not6}), $\lambda _{N}=\sum_{j=1}^{\infty
}\chi _{j}\left\vert c(j)+S_{N-1}(c_{j},Z)\right\vert ^{2}$. We use (\ref%
{NDCW7''}) with $\kappa=c(j)$ and $c$ replaced by $c_{j}$ (see (\ref{not5}))
and we obtain
\begin{eqnarray*}
{\mathbb{E}}_{V,\chi }(\lambda _{N}) &=&\sum_{j=1}^{\infty }\chi _{j}{%
\mathbb{E}}_{V,\chi }(\left\vert c(j)+S_{N-1}(c_{j},Z)\right\vert ^{2}) \geq
v^{N-1}(r)\sum_{j=1}^{\infty }\chi _{j}\sum_{\left\vert \alpha \right\vert
=N-1}c_{j}^{2}(\alpha )\chi ^{\alpha } \\
&=&v^{N-1}(r)\sum_{j=1}^{\infty }\chi _{j}\sum_{\left\vert \alpha
\right\vert =N-1}c^{2}(\alpha ,j)\chi ^{\alpha } =\frac{1}{N}%
v^{N-1}(r)\sum_{\left\vert \beta \right\vert =N}c^{2}(\beta )\chi ^{\beta }.
\end{eqnarray*}%
$\square $

We are now able to give our estimate:

\begin{lemma}
Let $c\in \mathcal{C}.$ For every $\eta >0$,
\begin{equation}
{\mathbb{P}}(\lambda _{N}\leq \eta )\leq C_{N}(r,\varepsilon )\Big(\Big(%
\frac{\eta}{|c|_{N}^{2}}\Big)^{1/N}+\exp \Big(-\frac{c_{N}(r,\varepsilon
)\left\vert c\right\vert _{N}^{2}}{\delta _{N}^{2}(c)}\Big)\Big)  \label{b1}
\end{equation}%
with $C_{N}(r,\varepsilon )$ denotes a constant of the type (\ref{C2}) and $%
c_{N}(r,\varepsilon )$ is given in (\ref{C1}).
\end{lemma}

\textbf{Proof.} We chose $J$ sufficiently large in order to have%
\begin{equation}
\sum_{\alpha \in \Gamma _{N}(J)}c^{2}(\alpha )\geq \frac{1}{2}\sum_{\alpha
\in \Gamma _{N}}c^{2}(\alpha )=\frac{1}{2}\left\vert c\right\vert _{N}^{2}.
\label{NDCW4a}
\end{equation}

We will use the Carbery--Wright inequality that we recall here (see Theorem
8 in \cite{[CW]}). Let $\mu $ be a probability law on ${\mathbb{R}}^{J}$
which is absolutely continuous with respect to the Lebesgue measure and has
a log-concave density. There exists a universal constant $K$ such that for
every polynomial $Q(x)$ of order $N$ and for every $\eta >0$ one has%
\begin{equation}
\mu (x:\left\vert Q(x)\right\vert \leq \eta )\leq KN(\eta /V_{\mu }(Q))^{1/N}
\label{NCW1}
\end{equation}%
with $V_{\mu }(Q)=(\int Q^{2}(x)d\mu (x))^{1/2}.$

We will use this result in the following framework. We recall that ${\mathbb{%
P}}_{V,\chi }$ is the conditional probability with respect to $\sigma
(V_{i},\chi _{i},i\in {\mathbb{N}})$ and we look to
\begin{equation*}
Q(U_{1},...,U_{J}):=\sum_{j=1}^{\infty }\chi _{j}\left\vert \partial
_{Z_{j}}S_{N}(c1_{\Gamma (J)},Z)\right\vert ^{2}=:\lambda _{N,J}
\end{equation*}%
as to a polynomial of order $N$ of $U_{1},...,U_{J}.$ It is easy to see that
the density of the law $\mu $ of $(U_{1},...,U_{J})$ (under ${\mathbb{P}}%
_{V,\chi })$\ is log-concave. So we are able to use (\ref{NCW1}). Using (\ref%
{NDC6})%
\begin{eqnarray*}
V_{\mu }(Q) &=&\Big(\int Q^{2}(x)d\mu (x)\Big)^{1/2}\geq \int \left\vert
Q(x)\right\vert d\mu (x)={\mathbb{E}}_{V,\chi }(\sum_{j=1}^{\infty }\chi
_{j}\left\vert \partial _{Z_{j}}S(c1_{\Gamma (J)},Z)\right\vert ^{2}) \\
&\geq &\frac{v^{N-1}(r)}{N}\sum_{\left\vert \beta \right\vert =N}c^{2}(\beta
)1_{\Gamma (J)}(\beta )\chi ^{\beta }.
\end{eqnarray*}%
We take now $\theta >0$ (to be chosen in a moment) and we use (\ref{NCW1})
in order to obtain%
\begin{eqnarray}
{\mathbb{P}}(\lambda_{N,J}\leq \eta)&=&{\mathbb{P}}(Q(U_{1},...,U_{J}) \leq
\eta )  \notag \\
&\leq& {\mathbb{P}}(V_{\mu }(Q)\leq \theta )+{\mathbb{E}}({\mathbb{P}}%
_{V,\chi }(Q(U_{1},...,U_{J})\leq \eta )1_{\{V_{\mu }(Q)\geq \theta \}})
\notag \\
&\leq &{\mathbb{P}}\Big(\sum_{\left\vert \beta \right\vert =N}c^{2}(\beta
)1_{\Gamma (J)}(\beta )\chi ^{\beta }\leq \frac{\theta N}{v^{N-1}(r)}\Big)%
+KN(\eta /\theta )^{1/N}.  \label{landa}
\end{eqnarray}
The first term in the above right hand side is estimated in Appendix \ref%
{hoef}: we apply Lemma \ref{H} with $x=\theta N/v^{N-1}(r)$ and with
the coefficients $\overline{c}_{J}(\alpha )=c(\alpha )1_{\Gamma
_{N}(J)}(\alpha ),$ so that $S_{N}(\overline{c}_{J}^{2},\chi )=$ $%
\sum_{|\beta| =N}c^{2}(\beta )1_{\Gamma (J)}(\beta )\chi ^{\beta }$. By (\ref%
{NDCW4a}) we have
\begin{equation*}
\frac{\left\vert c\right\vert _{N}^{2}}{2}\leq \left\Vert \overline{c}%
_{J}\right\Vert _{N}^{2}=\left\vert \overline{c}_{J}\right\vert _{N}^{2}\leq
\left\vert c\right\vert _{N}^{2}.
\end{equation*}

We recall that in Lemma \ref{H} we use $p=\varepsilon m(r)$ and that we need
(see (\ref{h2})) that
\begin{equation}
\theta =\frac{v^{N-1}(r)}{N}x\leq \frac{v^{N-1}(r)}{2N}\Big(\frac{p}{4}\Big)%
^{2N}\left\vert c\right\vert _{N}^{2}.  \label{theta}
\end{equation}%
We take $\theta $ equal to the quantity in the right hand side of the above
inequality so that
\begin{equation*}
x=\frac{\theta N}{v^{N-1}(r)}=\frac{1}{2}\Big(\frac{p}{4}\Big)%
^{2N}\left\vert c\right\vert _{N}^{2}.
\end{equation*}%
Then (\ref{h3}) gives%
\begin{equation*}
{\mathbb{P}}\Big(S_{N}(\overline{c}_{J}^{2},\chi )\leq\frac{\theta N}{%
v^{N-1}(r)}\Big)
\leq \frac{2e^{3}}{9}N\exp \Big(-\frac{1}{4}\Big(\frac{p}{4}\Big)%
^{4N}|c|_{N}^{4}\frac{1}{N\overline{\delta }_{N}^{2}(\overline{c}_{J})\|
\overline{c}_{J}\| _{N}^{2}}\Big).
\end{equation*}
%
%
%
Since $\overline{\delta }_{N}^{2}(\overline{c}_{J})\leq \delta _{N}^{2}(c)\ $%
and $\|\overline{c}_{J}\| _{N}^{2}\leq |c|_{N}^{2}$ we upper bound the above
term with
\begin{equation*}
\frac{2e^{3}}{9}N\exp \Big(-\frac{1}{4N}\Big(\frac{p}{4}\Big)^{4N}\frac{%
\left\vert c\right\vert _{N}^{2}}{\delta _{N}^{2}(c)}\Big).
\end{equation*}%
Inserting this in (\ref{landa}) we obtain

\begin{eqnarray*}
{\mathbb{P}}(\lambda _{N} &\leq &\eta )\leq {\mathbb{P}}(\lambda _{N,J}\leq
\eta ) \\
&\leq &\frac{2e^{3}}{9}N\exp \Big(-\frac{1}{4N}\Big(\frac{\varepsilon m(r)}{4%
}\Big)^{4N}\frac{|c| _{N}^{2}}{\delta _{N}^{2}(c)}\Big)+\frac{KN}{%
v(r)\varepsilon ^{2}m^{2}(r)| c| _{N}^{2/N}}\eta ^{1/N}.
\end{eqnarray*}%
and the proof is completed. $\square $

\subsection{Proof of the main results}\label{proofs}

Our basic lemma is the following:

\begin{lemma}
Let $Z,\overline{Z}\in \mathcal{L}((M_{p})_{p\geq 1},r,\varepsilon )$ and $c$%
,$\overline{c}\in \mathcal{C}.$\ We denote $S_{N}=S_{N}(c,Z)$ and $\overline{%
S}_{M}=S_{M}(\overline{c},\overline{Z}).$ Let $p_{\ast }$ be the universal
constant from (\ref{not14}). For every $k\in {\mathbb{N}}$ there exist a
constant $C_{N\vee M}(r,\varepsilon )$ as in (\ref{C2}) such that

\begin{equation}\label{b3}
\begin{array}{rl}
d_{0}(S_{N},\overline{S}_{M})
\leq &
\displaystyle
C_{N\vee M}(r,\varepsilon )%
(1+\Vert c\Vert _{N}+\Vert \overline{c}\Vert _{M})\Big(\frac{d_k^{\frac 1{k+1+kp_\ast N\vee M}}(S_N,\overline{S}_M)}{(|c|_N^{2/N}\wedge |\overline{c}|_M^{2/M})^{\frac {kp_\ast N\vee M}{k+1+kp_\ast N\vee M}}}\smallskip\\
&\displaystyle
+\exp \Big(-\frac{c_{N}(r,\varepsilon )\left\vert c\right\vert
_{N}^{2}}{\delta _{N}^{2}(c)}\Big)+\exp \Big(-\frac{c_{M}(r,\varepsilon )\left\vert
\overline{c}\right\vert _{M}^{2}}{\delta _{M}^{2}(\overline{c})}\Big)\Big),
\end{array}
\end{equation}
$c_N(\varepsilon, r)$ being given in (\ref{C1}).
\end{lemma}

\textbf{Proof.} We use Lemma \ref{d0dk} and in the estimate (\ref{not16}),
we replace ${\mathbb{P}}(\det \sigma _{S_{N}(c,Z)}<\eta )$ by the expression
from (\ref{b1}). So, we obtain%
\begin{equation*}
\begin{array}{rcl}
d_{0}(S_{N},\overline{S}_{M}) & \leq & \displaystyle C_{N\vee M}(r,\varepsilon )%
(1+\Vert c\Vert _{N}+\Vert \overline{c}\Vert _{M})\Big(\frac{1}{\eta ^{\frac{kp_{\ast}}{k+1}}}d_{k}^{\frac{1}{k+1}}(S_{N},\overline{S}%
_{M})+\smallskip \\
&  & +\frac 1{|c|_N^{2/N}\wedge |c|_M^{2/M}}\big(\eta^{1/N}+\eta^{1/M}\big)
+\exp (-\frac{c_{N}(r,\varepsilon )\left\vert c\right\vert
_{N}^{2}}{\delta _{N}^{2}(c)})+\exp (-\frac{c_{M}(r,\varepsilon )\left\vert
\overline{c}\right\vert _{M}^{2}}{\delta _{M}^{2}(\overline{c})})\Big).%
\end{array}%
\end{equation*}%
This holds true for every $\eta >0.$ We optimize over $\eta $ and we obtain (%
\ref{b3}). $\square $

\medskip

\textbf{Proof of Theorem \ref{CONVERGENCE}.} Let $S_{N}(c_{n},Z^{n}),n\in {%
\mathbb{N}}$, be the sequence considered in the statement of the theorem.
Since this sequence converges in law to $\mu ,$ it follows that it is a Cauchy sequence
in $d_{1}.$ And since $\overline{\delta }_{N}(c_{n})\rightarrow 0,$ and $%
\liminf_{n\rightarrow \infty }\left\vert c\right\vert _{N}^{2}>0$ the
inequality (\ref{b3}) says that the sequence is Cauchy in $d_{0}.$ It
follows that it converges to $\mu $ in $d_{0}.$ $\square $

\medskip

\textbf{Proof of Theorem \ref{INVARIANCE}.} By Theorem \ref{Smooth}, $%
d_{3}(S_{N}(c,Z),S_{N}(c,\overline{Z}))\leq C_{N}(r,\varepsilon )\overline{%
\delta }_{N}(c)$, so, using (\ref{b3}) with $k=3$ and $N=M$ we obtain (\ref{Int6}). $\square $

\medskip

\textbf{Proof of Theorem \ref{DISTANCES-2}.}
The hypotheses ensure that $\lim_{n}\overline{%
\delta }_{M}(c_{n})=0$, $\limsup_{n}\Vert c_{n}\Vert _{M}= \overline{C}_{M,X}  $,
$\liminf_{n}\Vert c_{n}\Vert _{M}= \underline{C}_{M,X}>0$
and $\lim_{n}d_{1}(S_{M}(c_{n},\overline{Z}),X )=0$. Notice that by Theorem \ref{CONVERGENCE}
we know that $\lim_{n}d_{0}(S_{M}(c_{n},\overline{Z}),X)=0.$
We write
\begin{align*}
&d_{0}(S_N(c,Z),X )
\leq d_{0}(S_N(c,Z),S_M(c_n,\overline{Z}))+d_{0}(S_M(c_n,\overline{Z}),X ) \\
&\quad \leq
C_{N\vee M}(r,\varepsilon )%
(1+\Vert c\Vert _{N}+\Vert c_n\Vert _{M})\Big(\frac{d_1^{\frac 1{2+p_\ast N\vee M}}(S_N(c,Z),S_M(c_n,\overline{Z}))}{(|c|_N^{2/N}\wedge |c_n|_M^{2/M})^{\frac {p_\ast N\vee M}{2+p_\ast N\vee M}}}\smallskip\\
&\quad
+\exp \Big(-\frac{c_{N}(r,\varepsilon )\left\vert c\right\vert
_{N}^{2}}{\delta _{N}^{2}(c)}\Big)+\exp \Big(-\frac{c_{M}(r,\varepsilon )\left\vert
c_n\right\vert _{M}^{2}}{\delta _{M}^{2}(c_n)}\Big)\Big)
+d_{0}(S_M(c_n,\overline{Z}),X )
\end{align*}%
the second inequality being (\ref{b3}).
Since $d_{1}(S_M(c_n,\overline{Z}),X )\rightarrow
0 $ then $d_{1}(S(c,Z),S_M(c_n,\overline{Z})\rightarrow$ $ d_{1}(S(c,Z),X ).$ We also have $\exp (-%
\frac{c_{M}(r,\varepsilon )\left\vert c_{n}\right\vert _{M}^{2}}{\delta
_{M}^{2}(c_{n})})+d_{0}(S_M(c_n,\overline{Z}),X)\rightarrow 0$ so (\ref{Int9-2}) is proved.
$\square $

\medskip

\textbf{Proof of Corollary \ref{DISTANCES}.}
Since $X \in \mathcal{A}%
_{M}((M_{p})_{p\geq 1},r,\varepsilon )$ we may find a sequence $%
Z^{(n)}=(Z_{k}^{(n)})_{k\in {\mathbb{N}}}\in \mathcal{L}((M_{p})_{p\geq
1},\varepsilon ,r)$ and a sequence $(c_{n})_n\subset \mathcal{C}$ that verifies the requests of Theorem \ref{DISTANCES-2}. So, the statement holds by repeating the proof of Theorem \ref{DISTANCES-2}.
$\square $

\section{Examples}

\subsection{Approximation with a chi-squared law} 

In \cite{[NP2009]}, Nourdin and Peccati give sufficient conditions in order
to estimate the Fortet-Mourier distance ($d_{1}$ in our notation) between a
multiple Wiener integral and a random variable with a centred Gamma
distribution. It is not clear if the Gamma distribution with fractional
coefficient is attainable in the sense of Definition \ref{Ateignable}, so we
are not able to use our results in the general case. But for an integer
parameter $\nu =2m$, the Gamma distribution coincides with the $\chi ^{2}$
distribution with $m$ degrees of freedom, and this law is clearly attainable
(just represent it as $\sum_{k=0}^{m-1}2\int_{k}^{k+1}W_{s}dW_{s}+m$ and
then use approximation with Riemann sums). So we restrict ourself to this
case. One looks to

\begin{equation*}
\Phi _{N}(c,Z)=\sum_{\alpha \in \Gamma _{N}}c(\alpha )Z^{\alpha }.
\end{equation*}%
If $Z_{k}$, $k\in {\mathbb{N}}$, are standard Gaussian random variables, then $%
\Phi _{N}(c,Z)$ is a multiple stochastic integral and in this case Nourdin
and Peccati in \cite{[NP2009]} have proved the following result. In order to
present it we have to introduce some notation. For $0\leq r\leq N$ and $%
\alpha ,\beta \in \Gamma _{N-r}$\ one denotes $c\otimes _{r}c(\alpha ,\beta
)=\sum_{\gamma \in \Gamma _{r}}c(\alpha ,\gamma )c(\beta ,\gamma )$ with the
convention that for $r=0$ we put $c\otimes _{0}c(\alpha ,\beta )=c(\alpha
)c(\beta )$ and for $r=N,$ $c\otimes _{N}c=\sum_{\gamma \in \Gamma
_{N}}c(\gamma )c(\gamma ).$ Notice that even if $c$ is symmetric, $c\otimes
_{r}c$ is not symmetric, so we introduce $c\widetilde{\otimes }_{r}c$ to be
the symmetrization of $c\otimes _{r}c.$ Finally, if $N$ is an even number,
we introduce
\begin{eqnarray*}
\kappa _{m,N}(c) &=&(m-N!\left\vert c\right\vert _{N}^{2})^{2}+4N!\left\vert
\theta _{N}\times c\widetilde{\otimes }_{N/2}c-c\right\vert _{N}^{2} \\
&&+N^{2}\sum_{\substack{ r\in \{1,...,N-1\}  \\ r\neq N/2}}%
(2N-2r)!(r-1)!^{2}\left(
\begin{array}{c}
N-1 \\
r-1%
\end{array}%
\right) ^{4}\left\vert c\otimes _{r}c\right\vert _{N}^{2}
\end{eqnarray*}%
with $\theta _{N}=\frac{1}{4}(N/2)!\left(
\begin{array}{c}
N \\
N/2%
\end{array}%
\right) .$ Combining Theorem 3.11 and Proposition 3.13 from \cite{[NP2009]} one
obtains the following:

\begin{theorem}
Let $N$ be an even integer and let $F(m)=\sum_{k=1}^{m}G_{k}^{2}-m$ with $%
G_{k}$ independent standard Gaussian random variables. Assume also that $%
Z_{k},k\in {\mathbb{N}}$ are independent standard Gaussian random variables.
Then
\begin{equation*}
d_{1}(\Phi _{N}(c,Z),F(m))\leq K_{1}(m)\kappa _{m,N}^{1/2}(c)
\end{equation*}%
with $K_{1}(m)=\max \{\sqrt{\pi /m},1/2m+1/2m^{2}\}.$
\end{theorem}

As an immediate consequence of Corollary \ref{DISTANCES} we obtain the
following result:

\begin{theorem}
\label{CHI2}Let $N$ be an even integer, $m\in {\mathbb{N}}_{\ast }$ and let $%
F(m)=\sum_{k=1}^{m}G_{k}^{2}-m$ with $G_{k}$ independent standard Gaussian
random variables. Assume also that $Z\in \mathcal{L}((M_{p})_{p\geq
1},r,\varepsilon ),$ and $c\in \mathcal{C}$. Then
$$
d_{0}(S_N(c,Z),X ) \leq
\displaystyle
C_{N\vee 2}(r,\varepsilon )(1+\| c\|_{N})
\Big(
\frac{\kappa _{m,N}^{\frac 1{4+2p_{\ast }N\vee M}}(c) }
{|c|_N^{\frac{2p_\ast N\vee 2}{N(2+p_{\ast }N\vee 2})}}
+\exp \Big(-\frac{c_{N}(r,\varepsilon )\left\vert c\right\vert _{N}^{2}}{%
\overline{\delta} _{N}^{2}(c)}\Big)\Big).
$$
\end{theorem}
\subsection{An example of quadratic CLT}\label{ex2}

An easy way to construct examples of invariance principles is to take a
double stochastic integral, to discretize it, and then to replace the
Brownian increments (renormalized) with some general random variables. So,
for example, starting with $\int_{0}^{1}\int_{0}^{t}f(t,s)dW_{s}W_{t}$ we
construct the approximation
\begin{equation*}
\sum_{0\leq i<j\leq 1}f\Big(\frac{i}{n},\frac{j}{n}\Big)\frac{\Delta _{i}}{\sqrt{n}}%
\frac{\Delta _{j}}{\sqrt{n}}
\end{equation*}%
with $\Delta _{i},i\in {\mathbb{N}}$ independent standard Gaussian random
variables. Then we replace $\Delta _{i}$ by some general $Z_{i}$ and we
obtain our invariance principle. Notice however that using this strategy
double sums give double integrals - so we remain in the same chaos. This is
true if $f$ is a square integrable function. In contrast, if we work with
some $f$ which is not square integrable then we may pass from a double sum
to a Gaussian limit (so to an element of the first chaos): a construction
phenomenon is at work. In this section we give an example which illustrates
this fact. We will study the convergence to normality of the following
stochastic series. We denote%
$$
S_{n}(Z)=\frac{1}{\sqrt{2n\ln n}}\sum_{i,j\geq 1}^n\frac{1}{\sqrt{\left\vert i-j\right\vert }}1_{\{i\neq j\}} Z_iZ_j.
$$
Notice that
\begin{equation}\label{quad1}
S_{n}(Z)=S_2(c_n,Z)\quad \mbox{with}\quad c_{n}(i,j)=\frac{%
1}{\sqrt{2n\ln n}}\frac{1}{\sqrt{\left\vert i-j\right\vert }}1_{\{i\neq j\}}1_{(i,j)\in\Gamma_2(n)}
\end{equation}

\begin{theorem}
\textbf{A}. Let $Z=(Z_{i})_{i\in {\mathbb{N}}}\in \mathcal{L}%
((M_{p})_{p},r,\varepsilon )$ and let $G=(G_{i})_{i\in {\mathbb{N}}}$ be a
sequence of standard normal random variables. Then%
\begin{equation}
d_{0}(S_{n}(Z),S_{n}(G))\leq \frac{C_{2}(r,\varepsilon )}{n^{1/(4+6p_{\ast
})}}  \label{quad8}
\end{equation}%
where $C_{2}(r,\varepsilon )=C(M_{p}r^{-1}\varepsilon ^{-1})^{q}$ with some
universal constants $C,p,q,$ and $p_{\ast }$ is the universal constant from (%
\ref{not14})

\textbf{B}. Let $W$ be a standard normal random variable. Then%
\begin{equation}
d_{0}(S_{n}(Z),W)\leq \frac{C_{2}(r,\varepsilon )}{\ln n}.  \label{quad10}
\end{equation}
\end{theorem}

\begin{remark}
Using the strategy mentioned in the beginning of this section we may easily
prove that, for $p<\frac{1}{2},$%
\begin{equation*}
\frac{1}{n^{1-p}}\sum_{i<j\leq n}\frac{1}{(j-i)^{p}}Z_{i}Z_{j}=\sum_{i<j\leq
n}\frac{1}{(\frac{j}{n}-\frac{i}{n})^{p}}\frac{Z_{i}}{\sqrt{n}}\frac{Z_{j}}{%
\sqrt{n}}\quad \overset{\mathcal{L}}{\longrightarrow }\quad
\int_{0}^{1}\int_{0}^{t}\frac{1}{(t-s)^{p}}dW_{s}dW_{t}.
\end{equation*}%
Notice that in this case we start with the function $f(t,s)=\left\vert
t-s\right\vert ^{-p}$ which is square integrable for $p<\frac{1}{2}.$ So,
with a soft singularity ($p<1/2$) we remain in the second chaos. But with a
strong singularity ($p=1/2$), a contraction phenomenon is at work and we
pass in the first chaos.
\end{remark}

\textbf{Proof.} \textbf{A.}
We apply Theorem \ref{INVARIANCE}. Here, $N=2$, $\|c_n\|_2=|c_n|_2$ and $\overline{\delta}_2(c_n)=\delta_2(c_n)$. So, by using (\ref{Int6}) we have
$$
d_0(S_n(Z),S_n(G))
\leq C_{3}(r,\varepsilon )(1+\|c_n\|_2)\Big(\frac{\overline{\delta }_{2}^{\frac{1}{4+6p_\ast }}(c_n)}{|c_n|_2^{\frac{6p_\ast}{4+6p_\ast }}}
+\exp \Big(-\frac{%
c_{2}(r,\varepsilon )\left\vert c_n\right\vert _{2}^{2}}{\overline{\delta }%
_{2}^{2}(c_n)}\Big)\Big),
$$
and (\ref{quad8}) immediately follows by using the estimates in (\ref{comp9}) e (\ref{comp9'}).

\smallskip

\textbf{B.} Let us prove (\ref{quad10}). We notice
that
\begin{equation*}
S_{n}(G)\overset{\mathcal{L}}{=}\int_{0}^{\infty
}\int_{0}^{t}f_{n}(t,s)dW_{s}dW_{t}=I_{2}(f_{n}),
\end{equation*}%
where $(W_{t})_{t}$ denotes a Brownian motion and%
\begin{equation*}
f_{n}(s,t)=\sum_{i,j=1}^{n}c_{n}(i,j)1_{(i,i+1]}(s)1_{(j,j+1]}(s).
\end{equation*}
Then we can use the results in \cite{[NP2009]} and we have that $%
d_{}(S_{n}(G),W)\leq C\sqrt{\kappa (f_{n})}$ where $\kappa (f_{n})$ is the
fourth cumulant of $I_{2}(f_{n}).$ And since $\kappa (f_{n})\leq
C\|f_{n}\otimes _{1}f_{n}\|_{L^2}^{2}=C\sum_{i,j}(c_{n}\otimes
_{1}c_{n})^{2}(i,j),$ (\ref{quad10}) is a consequence of (\ref{comp12}).
$\square $

\subsection{A variance-type estimator}

We denote
\begin{equation*}
X_{i}=\frac{1}{\sqrt{n}}\sum_{\substack{ j=1, \\ j\neq i}}^{n}\frac{1}{\sqrt{%
\left\vert i-j\right\vert }}Z_{j}
\end{equation*}%
and we study the asymptotic behavior of
\begin{equation*}
V_{n}(Z)=\sum_{i=1}^{n}(X_{i}^{2}-{\mathbb{E}}(X_{i}^{2})).
\end{equation*}%
The limit will be given by the double stochastic integral%
\begin{equation*}
I_{2}(\phi )=\int_{0}^{1}\int_{0}^{1}\phi (t,s)dW_{t}dW_{s}
\end{equation*}%
where the function $\phi $ is defined in (\ref{comp2}):%
\begin{equation*}
\phi (t,s)=\int_{0}^{1}\frac{du}{\sqrt{\left\vert (t-u)(s-u)\right\vert }}%
=\pi +2\ln \frac{\sqrt{1-t}+\sqrt{1-s}}{\left\vert \sqrt{t}-\sqrt{s}%
\right\vert }
\end{equation*}

\begin{proposition}
Let $Z=(Z_{i})_{i\in {\mathbb{N}}}\in \mathcal{L}((M_{p})_{p},r,\varepsilon )
$ Then%
\begin{equation}
d_{0}(V_{n}(Z),I_{2}(\phi ))\leq C_{2}(r,\varepsilon )(\frac{\ln ^{2}n}{n}%
)^{1/4(1+2p_{\ast })}.  \label{Var2}
\end{equation}
\end{proposition}

\textbf{Proof.} In this proof we refer several times to some computations and estimates which are developed in Appendix \ref{computations}.

\textbf{Step 1}. We denote%
\begin{eqnarray*}
a(i,j) &=&1_{\{i\neq j\}}\left\vert i-j\right\vert ^{-1/2}, \\
\overline{c}_{n}(i,j) &=&\frac{1}{n}(a\otimes _{1}a)(i,j):=\frac{1}{n}%
\sum_{k=1}^{n}a(i,k)a(j,k).
\end{eqnarray*}
We recall that in (\ref{comp10}),(\ref{comp10'}) and (\ref{comp10a}) one
proves that
\begin{equation}
\delta _{2}^2(\overline{c}_{n})\leq \frac{C\ln ^{2}n}{n},\quad 0<c_{\ast }\leq
\left\vert \overline{c}_{n}\right\vert ^{2}\leq C\quad \mbox{and}\quad
\sum_{k=1}^{n}\overline{c}_{n}^{2}(k,k)\leq \frac{C\ln ^{2}n}{n}.
\label{Ineg}
\end{equation}

We decompose%
\begin{equation*}
V_{n}(Z)=V_{n}^{\prime }(Z)+V_{n}^{\prime \prime }(Z)
\end{equation*}%
with%
\begin{eqnarray*}
V_{n}^{\prime }(Z) &=&\frac{2}{n}\sum_{j<j^{\prime }}a\otimes
_{1}a(j,j^{\prime })Z_{j}Z_{j^{\prime }}=\Phi_2(\overline{c}_n,Z)\quad \mbox{and} \\
V_{n}^{\prime \prime }(Z) &=&\frac{1}{n}\sum_{j=1}^{n}\Big(%
\sum_{i=1}^{n}a^{2}(i,j)\Big)(Z_{j}^{2}-1).
\end{eqnarray*}%
Since $V_{n}^{\prime \prime }(Z)$ contains terms of the form $Z_{j}^{2}$ we
may not use directly the results from the previous sections, and we are
obliged to develop a slight variant of them.

\textbf{Step 2}. By (\ref{C3})
\begin{equation*}
d_{3}(V_{n}^{\prime }(Z),V_{n}^{\prime }(G))\leq C\overline{\delta }_{2}(%
\overline{c}_{n})=C\delta_2(\overline{c}_n)\leq \frac{C\ln ^{2}n}{n}.
\end{equation*}%
And by the isometry property%
\begin{eqnarray*}
{\mathbb{E}}(\vert V_{n}^{\prime \prime }(Z)\vert ^{2}) &=&\frac{1%
}{n^{2}}\sum_{j=1}^{n}\Big(\sum_{i=1}^{n}a^{2}(i,j)\Big)^{2}{\mathbb{E}}%
((Z_{j}^{2}-1)^{2}) \\
&\leq &\frac{1}{n^{2}}\max_{j}(\E(Z_{j}^{4})-1)\sum_{j=1}^{n}\Big(%
\sum_{i=1}^{n}a^{2}(i,j)\Big)^{2} \\
&=&\frac{1}{n^{2}}\max_{j}(\E(Z_{j}^{4})-1)\sum_{j=1}^{n}\overline{c}%
_{n}^{2}(j,j) \\
&\leq &\frac{C\ln ^{2}n}{n}.
\end{eqnarray*}%
So%
\begin{equation}
d_{3}(V_{n}(Z),V_{n}(G))\leq \frac{C\ln ^{2}n}{\sqrt{n}}.  \label{Var3}
\end{equation}

\textbf{Step 3}. We will use the stochastic calculus of variations for $%
V_{n}(Z)$ so we have to estimate the Sobolev norms and the covariance
matrix. First
\begin{equation}
\left\Vert \left\vert V_{n}(Z)\right\vert \right\Vert _{q,p}\leq \left\Vert
\left\vert V_{n}^{\prime }(Z)\right\vert \right\Vert _{q,p}+\left\Vert
\left\vert V_{n}^{\prime \prime }(Z)\right\vert \right\Vert _{q,p}\leq
C_{2}(r,\varepsilon )\left\vert \overline{c}_{n}\right\vert _{2}^{2}.  \label{Var4}
\end{equation}%
This is because the estimate of $\left\Vert \left\vert V_{n}^{\prime }(Z)\right\vert
\right\Vert _{q,p}$ is already given in (\ref{not12}) and the estimate of $%
\left\Vert \left\vert V_{n}^{\prime \prime }(Z)\right\vert \right\Vert _{q,p}
$ is analogous (it suffices to follow the computations in Proposition 5.3 and 5.4 in \cite{[Numelin1]}), so we skip it.

We estimate now the covariance matrix (scalar in our case) defined in (\ref%
{not6}):%
\begin{equation*}
\lambda _{V_{n}(Z)}=\sum_{k=1}^{n}\chi _{k}\left\vert \partial
_{Z_{k}}V_{n}(Z)\right\vert ^{2}.
\end{equation*}%
We have
\begin{eqnarray*}
\partial _{Z_{k}}V_{n}(Z) &=&2\sum_{i=1}^{n}X_{i}\partial _{Z_{k}}X_{i}=%
\frac{2}{\sqrt{n}}\sum_{i=1}^{n}X_{i}a(i,k)\chi _{k} \\
&=&\frac{2}{n}\chi _{k}\sum_{i=1}^{n}a(i,k)\sum_{j=1}^{n}a(i,j)Z_{j} \\
&=&2\chi _{k}\sum_{j=1}^{n}Z_{j}(\frac{1}{n}\sum_{i=1}^{n}a(i,k)a(i,j))=2%
\chi _{k}\sum_{j=1}^{n}Z_{j}\overline{c}_{n}(j,k)
\end{eqnarray*}%
so that%
\begin{equation*}
\lambda _{V_{n}(Z)}=4\sum_{k=1}^{n}\chi _{k}\Big\vert \sum_{j=1}^{n}Z_{j}%
\overline{c}_{n}(j,k)\Big\vert ^{2}.
\end{equation*}%
This expression is strongly similar to $\lambda _{S_{2}(\overline{c}_{n},Z)}$
defined in (\ref{not6}), but there is one difference: we do not have the
property $\overline{c}_{n}(j,j)=0.$ So we have to eliminate the diagonal
terms. We define $c_{n}^{\prime }(i,j)=1_{\{i\neq j\}}\overline{c}_{n}(i,j)$
and we use the inequality $(a+b)^{2}\geq \frac{1}{2}a^{2}-b^{2}$ in order to
obtain%
\begin{eqnarray*}
\lambda _{V_{n}(Z)} &\geq &2\sum_{k=1}^{n}\chi _{k}\Big\vert
\sum_{j=1}^{n}Z_{j}c_{n}^{\prime }(j,k)\Big\vert ^{2}-4\sum_{k=1}^{n}\chi
_{k}\overline{c}_{n}^{2}(k,k)Z_{k}^{2} \\
&=&2\lambda _{S_{2}(c_{n}^{\prime },Z)}-4\sum_{k=1}^{n}\chi _{k}\overline{c}%
_{n}^{2}(k,k)Z_{k}^{2}.
\end{eqnarray*}%
Using (\ref{Ineg}), for $n$ sufficiently large we have%
\begin{equation*}
|c_{n}^{\prime }| _{2}^{2}\geq \frac{1}{2}|\overline{c}_{n}| _{2}^{2}-\sum_{k=1}^{n}\overline{c}%
_{n}^{2}(k,k)\geq \frac{c_{\ast }}{2}-\frac{C\ln ^{2}n}{n}\geq \frac{c_{\ast
}}{4}.
\end{equation*}%
Then, by (\ref{b1}) first and by (\ref{Ineg})  then, for every $\eta >0,$
\begin{align*}
{\mathbb{P}}(\lambda _{S_{2}(c_{n}^{\prime },Z)}\leq \eta )
&\leq
C_{2}(r,\varepsilon )\Big(\Big(\frac\eta{|c'_n|_2^2}\Big) ^{1/2}+\exp \Big(-\frac{c_{2}(r,\varepsilon
)\left\vert c_{n}^{\prime }\right\vert _{2}^{2}}{\delta
_{2}^{2}(c_{n}^{\prime })}\Big)\Big)\\
&\leq C_{2}(r,\varepsilon )(\eta ^{1/2}+\exp
(-c_{2}(r,\varepsilon )n)).
\end{align*}%
And again by (\ref{Ineg})
\begin{equation*}
{\mathbb{E}}\Big(\sum_{k=1}^{n}\chi _{k}c_{n}^{2}(k,k)Z_{k}^{2}\Big)\leq
\sum_{k=1}^{n}c_{n}^{2}(k,k)\leq \frac{C\ln ^{2}n}{n}
\end{equation*}%
so that%
\begin{eqnarray*}
{\mathbb{P}}(\lambda _{V_{n}(Z)} &\leq &\eta )\leq {\mathbb{P}}\Big(2\lambda
_{S_{2}(c_{n}^{\prime },Z)}\leq 2\eta \Big)+{\mathbb{P}}\Big(\sum_{k=1}^{n}\chi
_{k}c_{n}^{2}(k,k)Z_{k}^{2}\geq \eta \Big) \\
&\leq &C_{2}(r,\varepsilon )\Big(\eta ^{1/2}+\exp (-c_{2}(r,\varepsilon
)n)+\frac{C\ln ^{2}n}{\eta n}\Big)
\\
&\leq &C_{2}(r,\varepsilon )\Big(\eta ^{1/2}+\frac{C\ln ^{2}n}{\eta n}\Big).
\end{eqnarray*}

\textbf{Step 4}. We have all the ingredients in order that the regularization Lemma \ref{L3} holds for $V_n(Z)$ and $V_n(G)$ and we can prove for both of them an estimate as in (\ref{not16}).  By using it, we obtain, for $\eta<1$,
\begin{align*}
d_{0}(V_{n}(Z),V_{n}(G))
&\leq C\Big(\frac{1}{\eta ^{3p_{\ast }/4}}d_{3}^{\frac{1}{%
4}}(V_{n}(Z),V_{n}(G))
+{\mathbb{P}}(\lambda _{V_{n}(Z)} <\eta )+{\mathbb{P}}(\lambda
_{V_{n}(G)}<\eta )\Big) \\
&\leq C\Big(\frac{1}{\eta ^{3p_{\ast }/4}}\Big(\frac{\ln n}{\sqrt{n}}\Big)^{1/4}+\eta
^{1/2}+\frac{\ln ^{2}n}{\eta n}\Big)
\leq C\Big(\frac{1}{\eta ^{p_{\ast }}}\Big(\frac{\ln n}{\sqrt{n}}\Big)^{1/4}+\eta
^{1/2}\Big).
\end{align*}%
We optimize over $\eta<1 $ and we obtain%
\begin{equation*}
d_{0}(V_{n}(Z),V_{n}(G))\leq C\Big(\frac{\ln n}{\sqrt{n}}\Big)^{1/4(2p_{\ast }+1)}.
\end{equation*}

\textbf{Step 5}. Here we set $G_k=n(W_{k/n}-W_{(k-1)/n})$, where $W_t$ denotes the Brownian motion on which $I_2(\phi)$ is written. We estimate%
\begin{equation*}
\left\Vert V_{n}^{\prime }(G)-I_{2}(\phi )\right\Vert
_{2}^{2}=\int_{0}^{1}\int_{0}^{1}\left\vert \psi _{n}(x,y)-\phi
(x,y)\right\vert ^{2}dxdy,
\end{equation*}%
where
\begin{equation*}
\psi _{n}(x,y)=a\otimes _{1}a(i,j)\quad \mbox{for}\quad x\in I_{i},y\in I_{j}.
\end{equation*}%
By (\ref{comp3})%
\begin{equation*}
\left\vert \psi _{n}(x,y)-\phi (x,y)\right\vert \leq c(\frac{1}{\sqrt{n}%
\sqrt{\left\vert x-y\right\vert }}+\frac{1}{n(x+y)})
\end{equation*}%
so that
\begin{equation*}
\left\Vert V_{n}^{\prime }(G)-I_{2}(\phi )\right\Vert _{2}^{2}\leq \frac{C}{n%
}.
\end{equation*}%
Since $\lim_{n}\left\Vert V_{n}^{\prime \prime }(G)\right\Vert _{2}=0$ we
conclude that $\lim_{n}V_{n}(G)=I_{2}(\phi )$ in $L^{2}.$

Let $m\geq n.$ Using exactly the same argument as above we obtain, as $\eta<1$,%
\begin{align*}
d_{0}(V_{n}(G),V_{m}(G))
&\leq C\Big(\frac{1}{\eta ^{p_{\ast }}}d_{1}^{\frac{1}{%
2}}(V_{n}(G),V_{m}(G))
+{\mathbb{P}}(\lambda _{V_{n}(G)} <\eta )+{\mathbb{P}}(\lambda
_{V_{m}(G)}<\eta )\Big) \\
&\leq C\Big(\frac{1}{\eta ^{p_{\ast }}}\frac{1}{n^{1/2}}+\eta ^{1/2}+\frac{\ln
^{2}n}{\eta n}\Big)
\leq C\Big(\frac{1}{\eta ^{p_{\ast }}}\frac{1}{n^{1/2}}+\eta ^{1/2}\Big).
\end{align*}%
We optimize for $\eta<1 $ in order to obtain%
\begin{equation*}
d_{0}(V_{n}(G),V_{m}(G))\leq \frac{C}{n^{1/2(2p_{\ast }+1)}}.
\end{equation*}%
So $V_{n}(G),n\in \N$ is a  Cauchy sequence  in $d_{0}$ and consequently
converges to some limit which has to be $I_{2}(\phi ).$ And the estimate of
the error is the one given above. $\square $

\appendix

\section{An iterated Hoeffding's inequality}

\label{hoef}

In this section we estimate ${\mathbb{P}}(S_{N}(c^{2},\chi )\leq x)$
with
\begin{equation*}
S_{N}(c^{2},\chi )=\sum_{m=1}^{m_{0}}\sum_{\alpha \in \Gamma
_{m}}c^{2}(\alpha )\chi ^{\alpha }.
\end{equation*}%
Essentially this amounts to an iterated application of Hoeffding's
inequality. In order to implement this strategy we will use an extension of
Hoeffding's inequality to martingales, due to Benktus \cite{[Be]}. We recall that$%
\left\Vert c\right\Vert _{N}^{2}=\sum_{m=1}^{N}\left\vert c\right\vert
_{m}^{2}=\sum_{1\leq \left\vert \alpha \right\vert \leq N}c^{2}(\alpha )$
and $\overline{\delta }^{2}(c)$ is defined in (\ref{Int3}).

\begin{lemma}
\label{H}Let $p={\mathbb{P}}(\chi _{j}=1)=\varepsilon m(r).$ If
\begin{equation}
x\leq (\frac{p}{4})^{2N}\left\Vert c\right\Vert _{N}^{2}  \label{h2}
\end{equation}%
Then%
\begin{equation}
{\mathbb{P}}(S_{N}(c^{2},\chi )\leq x)\leq \frac{2e^{3}}{9}N\exp (-\frac{%
x^{2}}{N\overline{\delta }^{2}(c)\left\Vert c\right\Vert _{N}^{2}}).
\label{h3}
\end{equation}
\end{lemma}

\textbf{Proof.} We proceed by recurrence on $N.$ If $N=1$ we have
\begin{eqnarray*}
{\mathbb{P}}(S_{N}(c^{2},\chi ) &\leq &x)={\mathbb{P}}(\sum_{j}c^{2}(j)\chi
_{j}\leq x) \\
&\leq &{\mathbb{P}}(p\sum_{j}c^{2}(j)\leq 2x)+{\mathbb{P}}%
(\sum_{j}c^{2}(j)(p-\chi _{j})\geq x).
\end{eqnarray*}%
Since
\begin{equation*}
\sum_{j}c^{2}(j)=\left\Vert c\right\Vert _{1}^{2}\geq (\frac{4}{p})^{2}x>%
\frac{2x}{p}
\end{equation*}%
the first term is zero (here comes on the hypothesis (\ref{h2})). And by
Hoeffding's inequality%
\begin{equation*}
{\mathbb{P}}(\sum_{j}c^{2}(j)(p-\chi _{j})\geq x)\leq \exp (-\frac{2x^{2}}{%
\sum_{j}c^{4}(j)}).
\end{equation*}%
Since
\begin{equation*}
\sum_{j}c^{4}(j)\leq \max_{j}c^{2}(j)\times \sum_{j}c^{2}(j)=\overline{%
\delta }_{1}^{2}(c)\left\Vert c\right\Vert _{1}^{2}
\end{equation*}%
our inequality is verified.

Suppose now that (\ref{h3}) holds for $N-1$ and let us prove it for $N.$ We
recall that $\Gamma _{m}(j)=\{\alpha =(\alpha _{1},...,\alpha _{m}):\alpha
_{i}\leq j\}$ and we denote $\Gamma _{m}^{o}(j)=\{\alpha \in \Gamma
_{m}(j):\alpha _{1}<\alpha _{2}<....<\alpha _{m}\}.$ We also set $\Gamma_m^o=\{\alpha \in \Gamma_m\,:\,\alpha _{1}<\alpha _{2}<....<\alpha _{m}\}$. We write%
\begin{eqnarray*}
S_{N}(c^{2},\chi ) &=&\sum_{m=1}^{N}m!\sum_{\alpha \in \Gamma
_{m}^{o}}c^{2}(\alpha )\chi ^{\alpha } \\
&=&\sum_{j=1}^{\infty }c^{2}(j)\chi _{j}+\sum_{m=2}^{N}m!\sum_{j=1}^{\infty
}\chi _{j}\sum_{\alpha \in \Gamma _{m-1}^{o}(j-1)}c^{2}(\alpha ,j)\chi
^{\alpha } \\
&=&\sum_{j=1}^{\infty }\chi _{j}(c^{2}(j)+H_{j}) \\
&=&A+pB
\end{eqnarray*}%
with%
\begin{equation*}
H_{j}=\sum_{m=2}^{N}m!\sum_{\alpha \in \Gamma _{m-1}^{o}(j-1)}c^{2}(\alpha
,j)\chi ^{\alpha }.
\end{equation*}%
and%
\begin{equation*}
A=\sum_{j=1}^{\infty }(\chi _{j}-p)(c^{2}(j)+H_{j}),\qquad
B=\sum_{j=1}^{\infty }(c^{2}(j)+H_{j}).
\end{equation*}%
We take $x$ which satisfies (\ref{h2}) and we write%
\begin{equation*}
{\mathbb{P}}(S_{N}(c^{2},\chi )\leq x)\leq {\mathbb{P}}(B\leq 2x/p)+{\mathbb{%
P}}(-A\geq x)=:b+a.
\end{equation*}

Let us estimate $b.$ For $\alpha =(\alpha _{1},...,\alpha _{m})$ we denote $%
\overline{\alpha }=\max_{j=1,...,m}\alpha _{j}$ and
\begin{equation*}
\overline{c}^{2}(\alpha )=m\sum_{j>\overline{\alpha }}c^{2}(\alpha ,j)
\end{equation*}%
and we write%
\begin{eqnarray*}
\sum_{j=1}^{\infty }H_{j} &=&\sum_{m=2}^{N}m!\sum_{\alpha \in \Gamma
_{m-1}^{o}}\sum_{j>\overline{\alpha }}^{\infty }c^{2}(\alpha ,j)\chi
^{\alpha }=\sum_{m=2}^{N}\sum_{\alpha \in \Gamma _{m-1}}\sum_{j>\overline{%
\alpha }}^{\infty }mc^{2}(\alpha ,j)\chi ^{\alpha } \\
&=&\sum_{m=2}^{N}\sum_{\alpha \in \Gamma _{m-1}}\overline{c}^{2}(\alpha
)\chi ^{\alpha }=S_{N-1}(\overline{c}^{2},\chi ).
\end{eqnarray*}%
It follows that
\begin{equation*}
B=\sum_{j=1}^{\infty }c^{2}(j)+S_{N-1}(\overline{c}^{2},\chi ).
\end{equation*}

\textbf{Case 1}. We suppose that
\begin{equation}
\sum_{j=1}^{\infty }c^{2}(j)\geq \frac{1}{2}\left\Vert c\right\Vert _{N}^{2}.
\label{h5}
\end{equation}%
By (\ref{h2})%
\begin{equation*}
\frac{2}{p}x\leq (\frac{p}{2})^{2N-1}\left\Vert c\right\Vert _{N}^{2}<\frac{1%
}{2}\left\Vert c\right\Vert _{N}^{2}\leq \sum_{j=1}^{\infty }c^{2}(j)
\end{equation*}%
so that
\begin{equation*}
b={\mathbb{P}}(\sum_{j=1}^{\infty }c^{2}(j)+S_{N-1}(\overline{c}^{2},\chi
)\leq \frac{2}{p}x)=0.
\end{equation*}

\textbf{Case 2}. We suppose that
\begin{equation}
\sum_{j=1}^{\infty }c^{2}(j)<\frac{1}{2}\left\Vert c\right\Vert _{N}^{2}.
\label{h6}
\end{equation}%
Then ignore $\sum_{j=1}^{\infty }c^{2}(j)$ and we write%
\begin{equation*}
b\leq {\mathbb{P}}(S_{N-1}(\overline{c}^{2},\chi )\leq \frac{2}{p}x).
\end{equation*}%
We will use the recurrence hypothesis. Before doing this, we verify that
\begin{equation}
\overline{\delta }_{N-1}^{2}(\overline{c})\leq \overline{\delta }%
_{N}^{2}(c)\qquad \mbox{and}\qquad \frac{1}{4}\left\Vert c\right\Vert _{N}^{2}\leq
\left\Vert \overline{c}\right\Vert _{N-1}^{2}\leq \left\Vert c\right\Vert
_{N}^{2}.  \label{h7}
\end{equation}

Let $m\geq 2.$ We have%
\begin{eqnarray*}
\delta _{m-1}^{2}(\overline{c}) &=&\max_{j}\sum_{\alpha \in \Gamma _{m-1}}%
\overline{c}^{2}(\alpha ,j)=\max_{j}m\sum_{\alpha \in \Gamma
_{m-1}}\sum_{i>j\vee \overline{\alpha }}c^{2}(\alpha ,j,i) \\
&\leq &\max_{j}\sum_{\beta \in \Gamma _{m}}c^{2}(\beta ,j)=\delta
_{m}^{2}(c).
\end{eqnarray*}%
Summing over $m$ we obtain $\overline{\delta }_{N-1}^{2}(\overline{c})\leq
\overline{\delta }_{N}^{2}(c).$

We write now%
\begin{eqnarray*}
\left\Vert \overline{c}\right\Vert _{N-1}^{2}
&=&\sum_{m=1}^{N-1}m!\sum_{\alpha \in \Gamma _{m}^{o}}\overline{c}%
^{2}(\alpha )=\sum_{m=1}^{N-1}m!\sum_{\alpha \in \Gamma
_{m}^{o}}m\sum_{i>\alpha _{m}}c^{2}(\alpha ,i) \\
&=&\sum_{m=1}^{N-1}m!m\sum_{\beta \in \Gamma _{m+1}^{o}}c^{2}(\beta
)=\sum_{m=1}^{N-1}\frac{m}{m+1}\sum_{\beta \in \Gamma _{m+1}}c^{2}(\beta
)\leq \left\Vert c\right\Vert _{N}^{2}.
\end{eqnarray*}%
And, since $\frac{m}{m+1}\geq \frac{1}{2},$ we use (\ref{h6}) and we obtain
\begin{equation*}
\left\Vert \overline{c}\right\Vert _{N-1}^{2}\geq \frac{1}{2}%
\sum_{m=1}^{N-1}\sum_{\beta \in \Gamma _{m+1}}c^{2}(\beta )=\frac{1}{2}%
(\left\Vert c\right\Vert _{N}^{2}-\sum_{j=1}^{\infty }c^{2}(j))\geq \frac{1}{%
4}\left\Vert c\right\Vert _{N}^{2}
\end{equation*}%
so (\ref{h7}) is proved.

We have to verify that $\overline{x}=\frac{2}{p}x$ verifies (\ref{h2}).
Using (\ref{h2}) for $x$ and (\ref{h7}) we obtain%
\begin{equation*}
\overline{x}\leq \frac{4}{p}x\leq (\frac{p}{4})^{2N-1}\left\Vert
c\right\Vert _{N}^{2}\leq (\frac{p}{4})^{2(N-1)}\frac{1}{4}\left\Vert
c\right\Vert _{N}^{2}\leq (\frac{p}{4})^{2(N-1)}\left\Vert \overline{c}%
\right\Vert _{N-1}^{2}.
\end{equation*}%
Now we may use (\ref{h3}) and (\ref{h7}) and we obtain (notice that $%
x^{2}\leq \overline{x}^{2})$%
\begin{eqnarray*}
{\mathbb{P}}(S_{N-1}(\overline{c}^{2},\chi ) &\leq &\overline{x})\leq \frac{%
2e^{3}}{9}(N-1)\exp (-\frac{\overline{x}^{2}}{(N-1)\overline{\delta }%
_{N-1}^{2}(\overline{c})\left\Vert \overline{c}\right\Vert _{N-1}^{2}}) \\
&\leq &\frac{2e^{3}}{9}(N-1)\exp (-\frac{x^{2}}{N\overline{\delta }%
_{N}^{2}(c)\left\Vert c\right\Vert _{N}^{2}}).
\end{eqnarray*}%
We conclude that in both Case 1 and Case 2 we have%
\begin{equation}
b\leq \frac{2e^{3}}{9}(N-1)\exp (-\frac{x^{2}}{\varepsilon
_{N}^{2}(c)\left\Vert c\right\Vert _{N}^{2}}).  \label{h8}
\end{equation}

We estimate now $a.$ We denote
\begin{equation*}
h_{j}=c^{2}(j)+\sum_{m=2}^{N}m!\sum_{\alpha \in \Gamma
_{m-1}^{o}(j-1)}c^{2}(\alpha ,j).
\end{equation*}%
Since $0\leq \chi ^{\alpha }\leq 1$ we have%
\begin{equation*}
0\leq c^{2}(j)+H_{j}\leq h_{j}.
\end{equation*}%
Notice that%
\begin{equation*}
h_{j}=c^{2}(j)+\sum_{m=2}^{N}m\sum_{\alpha \in \Gamma
_{m-1}(j-1)}c^{2}(\alpha ,j)\leq N\overline{\delta }_{N}^{2}(c)
\end{equation*}%
and%
\begin{eqnarray*}
\sum_{j=1}^{\infty }h_{j} &=&\sum_{j=1}^{\infty
}c^{2}(j)+\sum_{m=2}^{N}m\sum_{j=1}^{\infty }\sum_{\alpha \in \Gamma
_{m-1}(j-1)}c^{2}(\alpha ,j) \\
&=&\sum_{j=1}^{\infty }c^{2}(j)+\sum_{m=2}^{N}\sum_{\beta \in \Gamma
_{m}}c^{2}(\beta )=\left\Vert c\right\Vert _{N}^{2}.
\end{eqnarray*}%
In particular%
\begin{equation*}
\sum_{j=1}^{\infty }h_{j}^{2}\leq N\overline{\delta }_{N}^{2}(c)\left\Vert
c\right\Vert _{N}^{2}.
\end{equation*}

We use now Corollary 1.4 pg 1654 in Bentkus [Be] which asserts the
following: if $M_{k},k\in {\mathbb{N}}$ is a martingale such that $%
\left\vert M_{k}-M_{k-1}\right\vert \leq h_{k}$ almost surely, then, for
every $n\in {\mathbb{N}},$
\begin{equation*}
{\mathbb{P}}(M_{n}\geq x)\leq \frac{2e^{3}}{9}\exp (-\frac{x^{2}}{%
\sum_{j=1}^{n}h_{j}^{2}}).
\end{equation*}%
In our case this gives%
\begin{equation*}
a={\mathbb{P}}(\sum_{j=1}^{\infty }(p-\chi _{j})(c^{2}(j)+H_{j})\geq x)\leq
\frac{2e^{3}}{9}\exp (-\frac{x^{2}}{N\overline{\delta }_{N}^{2}(c)\left\Vert
c\right\Vert _{N}^{2}}).
\end{equation*}%
This, together with (\ref{h8}) yields%
\begin{equation*}
a+b\leq \frac{2e^{3}}{9}N\exp (-\frac{x^{2}}{N\overline{\delta }%
_{N}^{2}(c)\left\Vert c\right\Vert _{N}^{2}}).
\end{equation*}%
$\square $

\section{Computations around an integral}\label{computations}

In this section we compute the following integral:
\begin{equation}
\phi (x,y)=\int_{0}^{1}\theta _{x,y}(z)dz\quad with\quad \theta _{x,y}(z)=%
\frac{1}{\sqrt{\left\vert x-z\right\vert \left\vert y-z\right\vert }}.
\label{comp1}
\end{equation}%
We also discuss the approximation with Riemann sums. We fix $n\in {\mathbb{N}%
}_{\ast }$ and we denote $I_{i}=[\frac{i}{n},\frac{i+1}{n})$ and $x_{i}=%
\frac{i}{n}.$

\begin{lemma}\label{lemmaB1}
For $0<x<y<1$, it holds
\begin{equation}
\phi (x,y)=\pi +2\ln \frac{\sqrt{1-x}+\sqrt{1-y}}{\left\vert \sqrt{x}-\sqrt{y%
}\right\vert }.  \label{comp2}
\end{equation}%
Moreover, if $x\in I_{i}$ and $y\in I_{j}$ with $i<j$ then%
\begin{equation}
\Big\vert \phi (x,y)-\frac{1}{n}\sum_{\substack{ k=1  \\ k\neq i,k\neq j}}%
^{n}\theta _{x_{i},x_{j}}(x_{k})\Big\vert \leq \frac{16\sqrt{2}}{\sqrt{n}}%
\frac{1}{\sqrt{y-x}}+\frac{8}{n(x+y)}.  \label{comp3}
\end{equation}
\end{lemma}

\textbf{Proof. Step 1}.
We consider the decomposition%
$$
(z-x)(z-y)=z^{2}-z(x+y)+xy=\left( z-\frac{x+y}{2}\right) ^{2}-\frac{(y-x)^{2}%
}{4}
$$
and we write
\begin{align*}
\phi(x,y)
=\int_{(0,x)\cup(y,1)}
\frac{1}{\sqrt{\big(z-\frac{x+y}2\big)^2-\big(\frac{y-x}2\big)^2}}\, dz
+\int_{(x,y)}
\frac{1}{\sqrt{\big(\frac{y-x}2\big)^2-\big(z-\frac{x+y}2\big)^2}}\, dz
\end{align*}
By using the change of variable $t=z-\frac{x+y}2$ and the fact that
\begin{eqnarray*}
\int \frac{dt}{\sqrt{t^{2}-a^{2}}} &=&\ln \left\vert t+\sqrt{t^{2}-a^{2}}%
\right\vert +C,\quad t^{2}>a^{2} \\
\int \frac{dt}{\sqrt{a^2-t^{2}}} &=&\arcsin \frac{t}{a}+C\quad t^{2}<a^{2},
\end{eqnarray*}%
straightforward computations give (\ref{comp2}).

\smallskip

\textbf{Step 2}.
We set
\begin{align*}
&I_{1}(a)
=\int_{a}^{x}\theta _{x,y}(z)dz
=\ln \frac{y-x}{2}-\ln \left\vert a-\frac{x+y}{2}+\sqrt{(a-x)(a-y)%
}\right\vert, \quad 0<a<x, \\
&I_{2}^{\prime }(a)
=\int_{x}^{a}\theta _{x,y}(z)dz
=\arcsin \frac{a-\frac{x+y}{2}}{\frac{y-x}{2}}+\frac{\pi }{2},\quad x<a<y, \\
&I_{2}^{\prime \prime }(a)
=\int_{a}^{y}\theta _{x,y}(z)dz
=\frac{\pi }{2}-\arcsin \frac{a-\frac{x+y}{2}}{\frac{y-x}{2}},\quad x<a<y,\\
&I_{3}(a)
=\int_{y}^{a}\theta _{x,y}(z)dz
=\ln \left\vert a-\frac{x+y}{2}+\sqrt{(a-x)(a-y)}\right\vert -\ln
\frac{y-x}{2},\quad y<a<1.
\end{align*}%
The above formulas in the last right hand sides follows by using the decomposition and the change of variable as in Step 1.

We first estimate $I_{i}(a)$ for $a$ close to $x$ or
to $y.$ First we notice that for $x-\frac{1}{n}<a<x<y$
\begin{equation}
I_{1}(a)=\int_{a}^{x}\frac{dz}{\sqrt{(x-z)(y-z)}}\leq \frac{1}{\sqrt{y-x}}%
\int_{a}^{x}\frac{dz}{\sqrt{x-z}}=\frac{2\sqrt{x-a}}{\sqrt{y-x}}\leq \frac{2%
}{\sqrt{n}}\frac{1}{\sqrt{y-x}}  \label{comp7}
\end{equation}%
and for $x<a<\frac{x+y}{2}\wedge (x+\frac{1}{n})$%
\begin{equation}
I_{2}^{\prime }(a)=\int_{x}^{a}\frac{dz}{\sqrt{(x-z)(y-z)}}\leq \frac{1}{%
\sqrt{y-\frac{x+y}{2}}}\int_{x}^{a}\frac{dz}{\sqrt{x-z}}=\frac{2\sqrt{2(x-a)}%
}{\sqrt{y-x}}\frac{1}{\sqrt{n}}\frac{2\sqrt{2}}{\sqrt{y-x}}.  \label{comp8}
\end{equation}%
Similar estimates hold for $I_{2}^{\prime \prime }(a)$ and for $I_{3}(a).$

We are now ready to prove (\ref{comp3}). We decompose%
\begin{equation*}
S=\frac{1}{n}\sum_{\substack{ k=1  \\ k\neq i,k\neq j}}^{n}\theta
_{x_{i},x_{j}}(x_{k})=S^{\prime }+S^{\prime \prime }+S^{\prime \prime \prime
}
\end{equation*}%
with%
\begin{equation*}
S^{\prime }=\frac{1}{n}\sum_{k=1}^{i-1}\theta _{x_{i},x_{j}}(x_{k}),\quad
S^{\prime \prime }=\frac{1}{n}\sum_{k=i+1}^{j-1}\theta
_{x_{i},x_{j}}(x_{k}),\quad S^{\prime \prime }=\frac{1}{n}%
\sum_{k=j+1}^{n}\theta _{x_{i},x_{j}}(x_{k}).
\end{equation*}%
And we also decompose
\begin{equation*}
I=\int_{0}^{1}\theta _{x,y}(z)dz=I^{\prime }+I^{\prime \prime }+I^{\prime
\prime \prime }
\end{equation*}%
with%
\begin{equation*}
I^{\prime }=\int_{0}^{x}\theta _{x,y}(z)dz,\quad I^{\prime \prime
}=\int_{x}^{y}\theta _{x,y}(z)dz\quad I^{\prime \prime \prime
}=\int_{y}^{1}\theta _{x,y}(z)dz.
\end{equation*}%
Let use estimate $I^{\prime }-S^{\prime }.$ We have $x_{i}\leq x<x_{i+1}$
and $x_{j}\leq y<x_{j+1}$\ so that
\begin{equation*}
x_{i}-x_{k}\leq x-x_{k}\leq x_{i}-x_{k-1},\quad x_{j}-x_{k}\leq y-x_{k}\leq
x_{j}-x_{k-1}
\end{equation*}%
so that%
\begin{equation*}
\frac{1}{\sqrt{\left\vert x_{i}-x_{k-1}\right\vert \left\vert
x_{j}-x_{k-1}\right\vert }}\leq \frac{1}{\sqrt{\left\vert x-x_{k}\right\vert
\left\vert y-x_{k}\right\vert }}\leq \frac{1}{\sqrt{\left\vert
x_{i}-x_{k}\right\vert \left\vert x_{j}-x_{k}\right\vert }}.
\end{equation*}%
Since $z\mapsto \theta _{x,y}(z)$ is increasing for $0<z<x$ we have%
\begin{equation*}
\frac{1}{n}\sum_{k=0}^{i-1}\frac{1}{\sqrt{\left\vert x-x_{k}\right\vert
\left\vert y-x_{k}\right\vert }}\leq I^{\prime }\leq \frac{1}{n}%
\sum_{k=1}^{i-1}\frac{1}{\sqrt{\left\vert x-x_{k}\right\vert \left\vert
y-x_{k}\right\vert }}+\int_{x_{i}}^{x}\theta _{x,y}(z)dz.
\end{equation*}%
Combining this with the previous inequality one gets%
\begin{eqnarray*}
\frac{1}{n}\sum_{k=0}^{i-2}\frac{1}{\sqrt{\left\vert x_{i}-x_{k}\right\vert
\left\vert x_{j}-x_{k}\right\vert }} &\leq &\frac{1}{n}\sum_{k=0}^{i-1}\frac{%
1}{\sqrt{\left\vert x-x_{k}\right\vert \left\vert y-x_{k}\right\vert }} \\
&\leq &I^{\prime }\leq \frac{1}{n}\sum_{k=1}^{i-1}\frac{1}{\sqrt{\left\vert
x_{i}-x_{k}\right\vert \left\vert x_{j}-x_{k}\right\vert }}%
+\int_{x_{i}}^{x}\theta _{x,y}(z)dz
\end{eqnarray*}%
One also has
\begin{equation*}
\frac{1}{n}\frac{1}{\sqrt{\left\vert x_{i}-x_{i-1}\right\vert \left\vert
x_{j}-x_{i-1}\right\vert }}\leq \frac{1}{\sqrt{n}}\frac{1}{\sqrt{y-x}}
\end{equation*}%
so that finally we obtain%
\begin{equation*}
\frac{1}{n}\sum_{k=1}^{i-1}\frac{1}{\sqrt{\left\vert x_{i}-x_{k}\right\vert
\left\vert x_{j}-x_{k}\right\vert }}-\frac{1}{\sqrt{n}}\frac{1}{\sqrt{y-x}}%
\leq I^{\prime }\leq \frac{1}{n}\sum_{k=1}^{i-1}\frac{1}{\sqrt{\left\vert
x_{i}-x_{k}\right\vert \left\vert x_{j}-x_{k}\right\vert }}%
+\int_{x_{i}}^{x}\theta _{x,y}(z)dz
\end{equation*}%
which, together with (\ref{comp7}), yields
\begin{equation*}
\left\vert I^{\prime }-S^{\prime }\right\vert \leq \frac{1}{\sqrt{n}}\frac{1%
}{\sqrt{y-x}}+\int_{x_{i}}^{x}\theta _{x,y}(z)dz\leq \frac{1}{\sqrt{n}}\frac{%
3}{\sqrt{y-x}}
\end{equation*}%
In a similar way one checks that
\begin{equation*}
\left\vert I^{\prime \prime \prime }-S^{\prime \prime \prime }\right\vert
\leq \frac{1}{\sqrt{n}}\frac{3}{\sqrt{y-x}}.
\end{equation*}%
In order to estimate $\left\vert I^{\prime \prime }-S^{\prime \prime
}\right\vert $ we note that $z\mapsto \theta _{x,y}(z)$ is increasing
for $x<z<\frac{x+y}{2}$ and decreasing for $\frac{x+y}{2}<z<y.$ So using
similar arguments we obtain, with $x_{l}\leq \frac{x+y}{2}<x_{l+1},$

\begin{eqnarray*}
\left\vert I^{\prime \prime }-S^{\prime \prime }\right\vert &\leq &\frac{4}{%
\sqrt{n}}\frac{1}{\sqrt{y-x}}+\int_{x}^{x_{i+1}}\theta
_{x,y}(z)dz+\int_{x_{j}}^{y}\theta _{x,y}(z)dz+\int_{x_{l}}^{x_{l+1}}\theta
_{x,y}(z)dz \\
&\leq &\frac{10\sqrt{2}}{\sqrt{n}}\frac{1}{\sqrt{y-x}}%
+\int_{x_{l}}^{x_{l+1}}\theta _{x,y}(z)dz.
\end{eqnarray*}%
It is easy to check that, if $\frac{1}{n}\leq \frac{y-x}{4}$
\begin{equation*}
\int_{\frac{x+y}{2}-\frac{1}{n}}^{\frac{x+y}{2}+\frac{1}{n}}\theta
_{x,y}(z)dz\leq \frac{1}{n}\times \frac{8}{y+x}
\end{equation*}%
And if $\frac{1}{n}>\frac{y-x}{4}$ then $I^{\prime \prime }$ does not
appear, so the above integral does not exists. So%
\begin{equation*}
\left\vert I^{\prime \prime }-S^{\prime \prime }\right\vert \leq \frac{10%
\sqrt{2}}{\sqrt{n}}\frac{1}{\sqrt{y-x}}+\frac{8}{n(x+y)}.
\end{equation*}%
We put all these inequalities together and we obtain%
\begin{equation*}
\left\vert I-S\right\vert \leq \frac{16\sqrt{2}}{\sqrt{n}}\frac{1}{\sqrt{y-x}%
}+\frac{8}{n(x+y)}.
\end{equation*}%
$\square $

\medskip

We will use Lemma \ref{lemmaB1} in order to compute the following quantities
which appear in our calculus. We denote%
\begin{eqnarray*}
a(i,j) &=&1_{\{i\neq j\}}\frac{1}{\sqrt{\left\vert i-j\right\vert }},\quad
c_{n}(i,j)=\frac{1}{\sqrt{2n\ln n}}a(i,j) \\
\overline{c}_{n}(i,j) &=&\frac{1}{n}\sum_{k=1}^{n}a(i,k)a(j,k)=(2\ln
n)\times (c_{n}\otimes _{1}c_{n})(i,j).
\end{eqnarray*}%
We also recall that
\begin{equation*}
\delta _{2}^2(c_{n})=\max_{i}\sum_{j=1}^{n}c_{n}^{2}(i,j),\quad \left\vert
c_{n}\right\vert _{2}^{2}=\sum_{i,j=1}^{n}c_{n}^{2}(i,j).
\end{equation*}

\begin{lemma}
\textbf{A.} We have%
\begin{align}
\delta _{2}^2(c_{n}) &\leq\frac{2}{n}  \label{comp9} \\
1-\frac{1}{\ln n} \leq \left\vert c_{n}\right\vert _{2}^{2}&\leq 1+\frac{1}{%
\ln n}  \label{comp9'}
\end{align}%
\textbf{B}. Let%
\begin{equation}
c_{\ast }=\frac{1}{16}\int_{\{\left\vert x-y\right\vert \geq \frac{1}{4}%
\}}\phi ^{2}(x,y)dxdy>0.  \label{comp10''}
\end{equation}
Then, for $n\geq (\frac{32\sqrt{2}}{\pi })^{2}$ one has%
\begin{align}
\delta _{2}^2(\overline{c}_{n}) &\leq \frac{C\ln ^{2}n}{n},  \label{comp10} \\
c_{\ast } \leq \left\vert \overline{c}_{n}\right\vert _{2}^{2}&\leq C,
\label{comp10'} \\
\sum_{k=1}^{n}\overline{c}_{n}^{2}(k,k) &\leq \frac{C\ln ^{2}n}{n}.
\label{comp10a}
\end{align}%
where $C$ is a universal constant.
\end{lemma}

\textbf{Proof. } We will first check that%
\begin{equation}
\ln i+\ln (n-i)\leq \sum_{j=1}^{n}a^{2}(i,j)\leq 2+\ln i+\ln (n-i).
\label{comp11}
\end{equation}%
Let us denote $x_{i}=\frac{i}{n}$ so that%
\begin{equation*}
a(i,j)=\frac{1}{\sqrt{\left\vert i-j\right\vert }}=\frac{1}{\sqrt{n}}\times
\frac{1}{\sqrt{\left\vert x_{i}-x_{j}\right\vert }}
\end{equation*}%
and then%
\begin{eqnarray*}
\ln (n-i) &=&\int_{\frac{1}{n}+x_{i}}^{1}\frac{dy}{y-x_{i}}\leq
\sum_{j=i+1}^{n-1}\frac{1}{x_{j}-x_{i}}\times \frac{1}{n}%
=\sum_{j=i+1}^{n}a^{2}(i,j) \\
&\leq &1+\int_{\frac{1}{n}+x_{i}}^{1}\frac{dy}{y-x_{i}}=1+\ln (n-i).
\end{eqnarray*}%
and%
\begin{eqnarray*}
\ln i &=&\int_{0}^{x_{i}-1/n}\frac{dy}{y-x_{i}}\leq \sum_{j=0}^{i-1}\frac{1}{%
x_{j}-x_{i}}\times \frac{1}{n}=\sum_{j=0}^{i-1}a^{2}(i,j) \\
&\leq &1+\int_{0}^{x_{i}-1/n}\frac{dy}{y-x_{i}}=1+\ln i.
\end{eqnarray*}%
Summing these two inequalities we obtain (\ref{comp11}).

Since $\int_{0}^{1}\ln xdx=-1$ we have
\begin{equation*}
n(\ln n-1)\leq \sum_{i=1}^{n}\ln i=n(\ln n+\frac{1}{n}\sum_{i=1}^{n}\ln
\frac{i}{n})\leq n\ln n
\end{equation*}%
so that summing over $i$ in (\ref{comp11}) we obtain
\begin{equation*}
2n(\ln n-1)\leq \sum_{i=1}^{n}\sum_{j=1}^{n}a^{2}(i,j)\leq 2n+2n\ln n
\end{equation*}%
which gives (\ref{comp9'}). And by (\ref{comp11})
\begin{equation*}
\delta _{2}^2(c_{n})=\max_{i}\sum_{j=1}^{n}c_{n}^{2}(i,j)\leq \frac{2(1+\ln n)%
}{2n\ln n}\leq \frac{2}{n}.
\end{equation*}%
so (\ref{comp9}) is also proved.

We will nw check that
\begin{equation}
\sum_{i,j=1}^{n}(c_{n}\otimes _{1}c_{n})^{2}(i,j)\leq \frac{C}{\ln ^{2}n}.
\label{comp12}
\end{equation}%
We construct the function
\begin{equation*}
\psi _{n}(x,y)=(a\otimes _{1}a)(i,j)\quad for\quad x\in I_{i},y\in I_{j}
\end{equation*}%
so that%
\begin{equation*}
(c_{n}\otimes _{1}c_{n})^{2}(i,j)=\frac{1}{4n^{2}\ln ^{2}n}(a\otimes
_{1}a)^{2}(i,j)=\frac{1}{4\ln ^{2}n}\int_{I_{i}\times I_{j}}\psi
_{n}^{2}(x,y)dxdy.
\end{equation*}%
Recall the function $\phi $ defined (\ref{comp1}). Using (\ref{comp3})
\begin{align*}
\sum_{\left\vert i-j\right\vert \geq 2}\int_{I_{i}\times I_{j}}\psi
_{n}^{2}(x,y)dxdy
\leq &2\sum_{\left\vert i-j\right\vert \geq
2}\int_{I_{i}\times I_{j}}\phi ^{2}(x,y)dxdy\\
&+2\sum_{\left\vert
i-j\right\vert \geq 2}\int_{I_{i}\times I_{j}}\left\vert \psi _{n}(x,y)-\phi
(x,y)\right\vert ^{2}dxdy \\
\leq &2\int \phi ^{2}(x,y)dxdy+C\int_{\{\left\vert x-y\right\vert \geq
\frac{1}{n}\}}\frac{1}{n\left\vert y-x\right\vert }+\frac{1}{n^{2}(x+y)^{2}}%
dxdy\leq C.
\end{align*}%
So%
\begin{equation}
\sum_{\left\vert i-j\right\vert \geq 2}(c_{n}\otimes _{1}c_{n})^{2}(i,j)\leq
\frac{C}{\ln ^{2}n}.  \label{comp13}
\end{equation}%
And, for $j\in \{i-1,i,i+1\}$%
\begin{eqnarray*}
\sum_{i=1}^{n}(c_{n}\otimes _{1}c_{n})^{2}(i,j)
&=&\sum_{i=1}^{n}(\sum_{k=1}^{n}c_{n}(k,i)c_{n}(k,j))^{2} \\
&\leq
&\sum_{i=1}^{n}(\sum_{k=1}^{n}c_{n}^{2}(k,i))(\sum_{k=1}^{n}c_{n}^{2}(k,j))
\\
&\leq &\delta ^{2}(c_{n})\left\vert c_{n}\right\vert _{2}^{2}\leq \frac{4}{n}%
.
\end{eqnarray*}%
So (\ref{comp12}) is proved.

Let us now prove that
\begin{equation}
\frac{c_{\ast }}{4\ln ^{2}n}\leq \sum_{i,j=1}^{n}(c_{n}\otimes
_{1}c_{n})^{2}(i,j)  \label{comp14}
\end{equation}%
Using (\ref{comp3}) and (\ref{comp2})%
\begin{eqnarray*}
\psi _{n}(x,y) &\geq &\phi (x,y)-\left\vert \phi (x,y)-\psi
_{n}(x,y)\right\vert  \\
&\geq &\frac{1}{2}\phi (x,y)+\frac{\pi }{2}-\frac{16\sqrt{2}}{\sqrt{n}\sqrt{%
\left\vert x-y\right\vert }}-\frac{8}{n(x+y)}
\end{eqnarray*}%
Notice that, if $\left\vert x-y\right\vert \geq \frac{1}{4}$ then $x+y\geq
\frac{1}{4}$.$T$ Then, if $\sqrt{n}\geq \frac{256\sqrt{2}}{\pi }$ we have
\begin{equation*}
\frac{16\sqrt{2}}{\sqrt{n}\sqrt{\left\vert x-y\right\vert }}\leq \frac{64%
\sqrt{2}}{\sqrt{n}}\leq \frac{\pi }{4},\quad and\quad \frac{8}{n(x+y)}\leq
\frac{32}{n}\leq \frac{\pi }{4}
\end{equation*}%
so that $\psi _{n}(x,y)\geq \frac{1}{2}\phi (x,y).$ It follows that
\begin{eqnarray*}
\sum_{i,j=1}^{n}(c_{n}\otimes _{1}c_{n})^{2}(i,j) &=&\frac{1}{4\ln ^{2}n}%
\int \psi _{n}^{2}(x,y)dxdy \\
&\geq &\frac{1}{16\ln ^{2}n}\int_{\{\left\vert x-y\right\vert \geq \frac{1}{4%
}\}}\phi ^{2}(x,y)dxdy=\frac{c_{\ast }}{\ln ^{2}n}.
\end{eqnarray*}%
So (\ref{comp14}) is proved. And (\ref{comp10'}) follows from (\ref{comp13})
and (\ref{comp14}).

Let us prove (\ref{comp10}). We fix $i$ and we write
\begin{align*}
\sum_{j>i}(c_{n}\otimes _{1}c_{n})^{2}(i,j) &\leq
2\sum_{j>i}\int_{I_{i}\times I_{j}}\phi ^{2}(x,y)dxdy
+2\sum_{j>i}\int_{I_{i}\times I_{j}}\left\vert \psi _{n}(x,y)-\phi
(x,y)\right\vert ^{2}dxdy\\
&=:A+B.
\end{align*}%
We have%
$$
A \leq 2\int_{I_{i}\times I_{i+1}}\phi
^{2}(x,y)dxdy+2\sum_{j>i+1}\int_{I_{i}\times I_{j}}\phi ^{2}(x,y)dxdy
=:A^{\prime }+A^{\prime \prime }.
$$
Using (\ref{comp2}),%
\begin{align*}
A^{\prime }
&\leq \frac{C}{n^{2}}+C\int_{I_{i}\times I_{i+1}}\ln ^{2}\frac{1}{%
\left\vert x-y\right\vert }dxdy\leq \frac{C\ln ^{2}n}{n^{2}},\\
A^{\prime \prime }
&\leq C\int_{x_{i}}^{x_{i+1}}\int_{x_{i+1}+\frac{1}{n}%
}^{1}\phi ^{2}(x,y)dxdy\leq \frac{C\ln ^{2}n}{n}.
\end{align*}%
Using the estimate (\ref{comp3}) we get similar estimates for $B.$ And this
gives (\ref{comp10}).

We prove now (\ref{comp10a}). Using (\ref{comp9}) and (\ref{comp9'})
$$
\sum_{k=1}^{n}\overline{c}_{n}^{2}(k,k) =\sum_{k=1}^{n}\Big(\frac{1}{n}%
\sum_{i=1}^{n}a^{2}(i,k)\Big)^{2}
\leq \max_{i}\frac{1}{n}\sum_{i=1}^{n}a^{2}(i,k)\times \frac{1}{n}%
\sum_{k,i=1}^{n}a^{2}(i,k)\leq C\frac{\ln ^{2}n}{n}.
$$
$\square $


\addcontentsline{toc}{section}{References}

\begin{thebibliography}{99}
\bibitem{[BC-CLT]} V. Bally, L. Caramellino. Asymptotic development for the
CLT in total variation distance. \emph{Bernoulli} \textbf{22}, 2442--2485,
2016.

\bibitem{[BC-EJP]} V. Bally, L. Caramellino. On the distances between
probability density functions. \emph{Electronic Journal of Probability}
\textbf{19}, no. 110, 1--33, 2014.

\bibitem{[Numelin1]} V. Bally, L. Caramellino. An invariance principle for
stochastic series I. Gaussian limits. \texttt{ArXiv:1510.03616}, 2015.

\bibitem{[BCl]} V. Bally, E. Cl\'ement. Integration by parts formula and
applications to equations with jumps. \emph{Probab. Theory Related Fields},
\textbf{151}, 613--657, 2011.

\bibitem{[BGL]} D. Bakry, I. Gentil, M. Ledoux \emph{Analysis and Geometry
of Markov Diffusion Semigroups}. Springer, 2014.

\bibitem{[Be]}
V. Bentkus. On Hoeffding's inequalities. \emph{Ann. Probab.} \textbf{32},  1650--1673, 2004.

\bibitem{[BGJ]} K. Bichtler, J.-B. Gravereaux, J. Jacod. \emph{Malliavin
calculus for processes with jumps}. Gordon and Breach Science Publishers,
1987.

\bibitem{[CW]} A. Carbery, J. Wright. Distributional and $L^q$ norm
inequalities for polynomials over convex bodies in ${\mathbb{R}}^n$. \emph{%
Math. Research Lett.} \textbf{8}, 233--248, 2001.

\bibitem{[dJ1]} P. de Jong. A central limit theorem for generalized
quadratic forms. \emph{Probab. Th. Rel. Fields} \textbf{75}, 261--277, 1987.

\bibitem{[dJ2]} P. de Jong. A central limit theorem for generalized
multilinear forms. \emph{Journal of Multivariate Analysis} \textbf{34},
275--289, 1990.

\bibitem{Fish} R.A. Fisher. Moments and product moments of sampling
distributions. \emph{Proceedings of the London Mathematical Society} \textbf{%
2}, 199--238, 1929.

\bibitem{Hoe} W. Hoeffding. A class of statistics with asymptotically normal
distributions. \emph{Ann. Statistics} \textbf{19}, 293--325, 1948.

\bibitem{[IW]} N. Ikeda, S. Watanabe.  \emph{Stochastic Differential
Equations and Diffusion processes}. North-Holland Mathematical Library 24,
1989.

\bibitem{L} R. Latala. Estimates of moments and tails of Gaussian chaoses.
\emph{Ann. Probab.} \textbf{34}, 2315--2331, 2006.

\bibitem{lee} A.J. Lee. \emph{U-Statistics: Theory and Practice.} Marcel
Dekker, New York, 1990.

\bibitem{[MPy]} D. Malicet, G. Poly. Properties of convergence in Dirichlet
structures. \emph{J. Funct. Anal.} \textbf{264}, 2077--2096, 2013.

\bibitem{[MO'DO]} E. Mossel, R. O'Donnell, K. Oleszkiewicz. Noise stability
of functions with low influences: Variance and optimality. \emph{Ann. Math.}
\textbf{171}, 295--341, 2010.

\bibitem{[NN]} S. Noreddine, I. Nourdin. On the Gaussian approximation of
vector-valued multiple integrals. \emph{J. Multiv. Anal.} \textbf{102},
1008-1017, 2011.

\bibitem{[NP]} I. Nourdin, G. Peccati. \emph{Normal Approximations Using
Malliavin Calculus: from Stein's Method to Universality}. Cambridge Tracts
in Mathematics, 192, 2012.

\bibitem{[NP2009]} I. Nourdin, G. Peccati. Stein's method on Wiener chaos.
\emph{Probab. Theory Related Fields} \textbf{145}, 75--118, 2009.

\bibitem{[NPRein]} I. Nourdin, G. Peccati, G. Reinert. Invariance principles
for homogeneous sums: universality of Wiener chaos. \emph{Ann. Probab.}
\textbf{38}, 1947--1985, 2010.

\bibitem{[NPRev]} I. Nourdin, G. Peccati, A. R\'{e}veillac. Multivariate
normal approximation using Stein's method and Malliavin calculus. \emph{Ann.
Inst. H. Poincar\'{e} Probab. Statist.} \textbf{46}, no. 1, 45--58, 2010.

\bibitem{[NPy]} I. Nourdin, G.Poly. Convergence in total variation on Wiener
chaos. \emph{Stochastic Process. Appl.} \textbf{123}, 651--674, 2013.

\bibitem{[NPy2]} I. Nourdin, G.Poly. Convergence in law in the second
Wiener/Wigner chaos. Convergence in law in the second Wiener/Wigner chaos.
\emph{Elect. Comm. in Probab.} \textbf{17}, no. 36, 2012.

\bibitem{bib:[N]} D. Nualart. \emph{The Malliavin calculus and related
topics. Second Edition}. Springer-Verlag, 2006.

\bibitem{[NO]} D. Nualart, S. Ortiz-Latorre. Central limit theorem for
multiple stochastic integrals and Malliavin calculus. \emph{Stoch. Processes
Appl.} \textbf{118}, 614--628, 2008.

\bibitem{[PN]} D. Nualart, G. Peccati. Central limit theorems for sequences
of multiple stochastic integrals. \emph{Annals of Probability} \textbf{33},
177--193, 2005.

\bibitem{[PT]} G. Peccati, C.A. Tudor. Gaussian limits for vector-valued
multiple stochastic integrals. \emph{S\'{e}minaire de Probabilit\'{e}s
XXXVIII}, 247--262, 2004.

\bibitem{proh} Yu.V. Prohorov. A local theorem for densities. \emph{Doklady
Akad. Nauk SSSR (N.S.)} \textbf{83}, 797--800, 1952. 

\bibitem{Sev} B.A Sevastianov. The class of limit laws for distributions of
quadratic forms in normal variables. \emph{Theor. Probability Appl.} \textbf{%
6}, 368---372, 1961.
\end{thebibliography}
\end{document}